%% file: spectral_RQMC.tex
\documentclass{amsart}
\usepackage{graphicx}

\usepackage{orcidlink}
\usepackage{amsmath,amssymb}
\usepackage{bm}
\usepackage{multirow}

\usepackage{hyperref}
\usepackage[misc]{ifsym} 

\usepackage{color}
\usepackage{graphicx}
\usepackage{subcaption}

\usepackage{amsfonts}
\usepackage{mathtools}

\usepackage[mathscr]{eucal}

\usepackage{bbm}

\input{commands_blue_noise.tex}

\ifpdf
\DeclareGraphicsExtensions{.eps,.pdf,.png,.jpg}
\else
\DeclareGraphicsExtensions{.eps}
\fi

\begin{document}

\title{Randomized quasi-Monte Carlo and Owen's boundary growth condition: A spectral analysis}

\author[Y.~Liu]{Yang Liu}
\address[Y.~Liu]{CEMSE, King Abdullah University of Science and Technology, Thuwal, Saudi Arabia} \email[]{yang.liu.3@kaust.edu.sa}

\subjclass{65C05}

\begin{abstract}
{In this work, we analyze the convergence rate of randomized quasi-Monte Carlo (RQMC) methods under Owen's boundary growth condition [Owen, 2006] via spectral analysis. Specifically, we examine the RQMC estimator variance for the two commonly studied sequences: the lattice rule and the Sobol' sequence. We apply the Fourier transform and Walsh--Fourier transform, respectively. Assuming certain regularity conditions, our findings reveal that the asymptotic convergence rate of the RQMC estimator's variance closely aligns with the exponent specified in Owen's boundary growth condition for both sequence types. We also provide analysis for certain discontinuous integrands. }

\end{abstract}
\maketitle

\noindent \textbf{Keywords.} Randomized quasi-Monte Carlo; Variance; Lattice rule; Sobol' sequence; Fourier transform; Walsh–Fourier transform.

\section{Introduction}
\label{sec;introduction}

The quasi-Monte Carlo (QMC) method utilizes deterministic low-discrepancy sequences (LDS) to compute integrals. The Koksma--Hlawka inequality indicates that, for integrands whose variation in the Hardy--Krause sense is finite, the QMC integration error is bounded by the product of the Hardy--Krause variation and the star-discrepancy of the point set. However, the Koksma--Hlawka inequality fails to be informative for certain singular integrands in $L^2$. For functions with boundary unboundedness, \cite{owen2006halton} characterizes a boundary growth condition and establishes connections with the asymptotic QMC convergence rates. Meanwhile,~\cite{kuo2012quasi} find that problem-specific lattice rules can achieve dimension-independent $\mathcal{O}(N^{-1})$ worst-case convergence rate. Additionally, \cite{he2015convergence} study the convergence of randomly scrambled Sobol' sequence for integrands with discontinuities. 

The randomized QMC (RQMC) achieves unbiased mean estimators by randomizing the LDSs. The variance of the RQMC estimator has been explored in the previous literature. \cite{lecuyer2000variance} investigate the randomly-shifted lattice rule (RSLR) using the Fourier transform. \cite{owen1997monte} analyzes Sobol' sequences through Haar wavelets. An equivalent formulation using the Walsh--Fourier transform is revealed in Chapter 13 of~\cite{dick2010digital}. In this work, we conduct a spectral analysis of the RQMC estimator variance for both lattice rules and Sobol' sequences concerning certain types of singular integrands. Assuming these integrands exhibit boundary unboundedness, interior discontinuities, or both, we leverage Fourier and Walsh--Fourier transforms to investigate their spectral properties. 

This paper is organized as follows. Section~\ref{sec:background} provides background information and notations related to the RQMC. Section~\ref{sec:spectral_analysis_rqmc} analyzes the RSLR and scrambled Sobol' sequences through spectral analysis. Section~\ref{sec:discontinuous_integrand} further analyzes Sobol' sequences for discontinuous integrands. Section~\ref{sec:numex} presents numerical results. Finally, section~\ref{sec:conclusion} concludes the paper and discusses possible further directions. 

\section{Background}
\label{sec:background}
In this section, we introduce the lattice rule, Sobol' sequences, and Owen's boundary growth condition to characterize the integrand regularity. 
\subsection{Notations}
\label{sec:notations}
First, we want to introduce notations that we use throughout the paper. 
\begin{description}
	\item $A^{*} = \max_j A_j$: the maximum of the boundary growth condition parameters. 
	\item $[\bm{a}, \bm{b}]$: the hyperrectangle interval $\{ \bm{t} \in \mathbb{R}^s: \min({a}_j, {b}_j) \leq {t}_j \leq \max({a}_j, {b}_j) \}$. This notation is particularly convenient when considering the intervals of the form $[\bm{t}, \bm{t} + \frac{1}{2\bm{n}} \ \mathrm{mod} \ 1] = \{ \bm{\tau} \in \mathbb{R}^s: \min({t}_j, {t}_j + \frac{1}{2n_j} \ \mathrm{mod} \ 1) \leq {\tau}_j \leq \max({t}_j, {t}_j + \frac{1}{2n_j} \ \mathrm{mod} \ 1) \}$ for $\bm{n} \in \mathbb{N}^s$.
	\item ${}_{b} \mathrm{wal}_{\bm{\ell}}$: the Walsh basis with index $\bm{\ell} \in \mathbb{N}_0^s$ in base $b$. 
	\item $c_{\bm{n}}$: the Fourier coefficient of $f$, $\bm{n} \in \mathbb{Z}^s$.
	\item $\Delta(f;\cdot, \cdot)$: the alternating sum of $f$. 
	\item $\bar{f}(\bm{\ell})$: the Walsh--Fourier coefficient of $f$ with index $\bm{\ell} \in \mathbb{N}_0^s$. 
	\item $\phi:[0, 1]\to \mathbb{R}$ is given by $\phi(t) \coloneqq \min(t, 1-t)$. To be used in Owen's boundary growth condition and derivations. 
	\item $\Gamma_{\mathfrak{u}, \bm{\ell}}$, $\Gamma_{\bm{\ell}}$: the ``gain coefficient'' of the digital net, $\mathfrak{u} \subseteq 1:s$, $\bm{\ell} \in \mathbb{N}_0^s$. 
	\item $\lesssim$: the inequality holds up to a finite constant. $a \lesssim b$ means $a \leq C b$ for some constant $C < +\infty$, $a, b \in \mathbb{R}$. 
	\item $s \in \mathbb{N}$: the integration dimension. We denote $1:s = \{1, 2, \dotsc, s\}$. 
	\item $\bm{t}^{\mathfrak{u}} : \bm{\tau}^{-\mathfrak{u}} $: a concatenation of two vectors $\bm{t}$ and $\bm{\tau}$. The $j$-th component of this concatenated vector is given by $t_j$ if $j \in \mathfrak{u}$, and by $\tau_j$ otherwise. 
\end{description}
\subsection{Lattice rule}
In this work, we specifically discuss the rank-1 lattice rule. Given a generating vector $\bm{z} \in \mathbb{N}^s$, the lattice rule with $N$ quadrature points is given by
\begin{equation}
	\label{eq:lattice_rule}
	\left\{  i \frac{\bm{z}}{N} \ \textrm{mod} \ 1 : i = 0, 1, \dotsc, N-1 \right\}.
\end{equation}
The lattice rule can be extensible, meaning that when $N$ increases, new quadrature points can be added to the existing set without discarding the previously computed points. This property is particularly useful in practical numerical simulations. Extensible lattice rules are often addressed to as ``lattice sequences''. Analysis of lattice sequences can be found, for instance, in~\cite{hickernell2000extensible}. 

The choice of the generating vector $\bm{z}$ influences the integration error of the lattice rule. Approaches for optimizing the generating vector are discussed in works such as~\cite{kuo2003component, cools2006constructing} and more recently in~\cite{l2020tool}. 

\subsection{Sobol' sequence}
We first introduce the definition of a $(t, m, s)$-net in base $b$. A $(t, m, s)$-net in base $b$ is a set of $b^m$ points in $[0, 1)^s$ such that every $s$-dimensional elementary interval
	\begin{equation}
		\label{eq:elementary_interval}
		E_{\bm{\ell}, \bm{k}} = \prod_{j=1}^{s} \left[ \frac{{k}_j}{b^{{\ell}_j}},  \frac{{k}_j + 1}{b^{{\ell}_j}}\right),
	\end{equation}
	where $\bm{\ell} = ({\ell}_1, \dotsc, {\ell}_s) \in \mathbb{N}_0^s$, {$\abs{\bm{\ell}} := \sum_{j=1}^s {\ell}_j = m - t$} and $\bm{k} = ({k}_1, \dotsc, {k}_s) \in \mathbb{N}_0^s$ with ${k}_j < b^{{\ell}_j}$ for $j = 1, \dotsc, s$, contains exactly {$b^t$} points from the net. The smallest integer $t$ for which this property holds is referred to as the quality parameter; $s$ is the dimension and $b^m$ is the number of quadrature points. More properties of the $(t, m, s)$-nets can be found in~\cite{dick2010digital}. 
	
	The Sobol' sequence is an example of a $(t, m, s)$-net in base $b = 2$. The parameters of the Sobol' sequence can be optimized to improve the quadrature points quality. For instance, \cite{joe2008constructing} optimizes direction numbers of the Sobol' sequences to achieve better two-dimensional projections. 
	
	Regarding randomization techniques, \cite{owen1995randomly} proposes nested uniform scrambling to randomize the Sobol' sequence, which, however, requires substantial storage and computational cost. \cite{matouvsek1998thel2} proposes the random linear scramble as an alternative approach with reduced computational cost, which does not alter the mean squared $L^2$ discrepancy. A review of other randomizations can be found in~\cite{owen2003variance}. 
\subsection{Owen's boundary growth condition}
\cite{owen2006halton} proposes the boundary growth condition to characterize the unbounded integrand $f$ on $[0, 1]^s$, which is given as follows.
\begin{Assumption}[Owen's boundary growth condition~\cite{owen2006halton}]
	\label{assumption:owen_boundary_growth_condition}
	 Let $f \in L^2([0, 1]^s)$, and the derivative $\frac{\partial^{\mathfrak{u}} f}{\partial \bm{t}_{\mathfrak{u}}}$ exists on $(0, 1)^s$ for $\mathfrak{u} \subseteq 1:s$. Then the boundary growth condition is given by
	\begin{equation}
		\label{eq:boundary_growth_condition}
		\abs*{\frac{\partial^{\mathfrak{u}} f}{\partial \bm{t}_{\mathfrak{u} } }} \lesssim \prod_{j=1}^s \phi(t_j)^{-\mathbbm{1}_{j \in \mathfrak{u}} - A_j},
	\end{equation}
for all $\mathfrak{u} \subseteq \{1, \dotsc, s\}$, $\max_j A_j < 1/2$ and $\phi:[0, 1]\to \mathbb{R}$, $\phi(t) \coloneqq \min(t, 1-t)$. 
\end{Assumption}
\cite{owen2006halton} proves that for integrands satisfying the condition~\eqref{eq:boundary_growth_condition}, the $L^1$ convergence rate of RQMC using low-discrepancy sequences is $\mathcal{O}(N^{-1 + \epsilon + A^{*}})$ for $\epsilon > 0$ and $A^* = \max_j A_j$ as defined in Section~\ref{sec:notations}. \cite{he2022error} and~\cite{gobet2022mean} have derived the RQMC $L^2$ convergence rate $\mathcal{O}(n^{-1+\epsilon+A^{*}})$ for the scrambled Sobol' sequence. The former applies the Koksma--Hlawka inequality, while the latter establishes results based on the pairwise correlations on randomized QMC quadrature points. 

In this work, we investigate the $L^2$ convergence rates of the randomized lattice rule and Sobol' sequence. In the next section, we perform a spectral analysis of these two low-discrepancy sequences under the boundary growth condition assumptions.

\section{Spectral Analysis of RQMC and Owen's boundary growth condition}
\label{sec:spectral_analysis_rqmc}
In this section, we explore the variance of the RSLR and scrambled Sobol' sequence using the Fourier and Walsh--Fourier transformations, respectively. 

\subsection{Lattice Rule}
We consider the integrand to be periodic. Let $f : \mathbb{T}^s = [0, 1]^s \to \mathbb{R}$ and $f \in L^2(\mathbb{T}^s)$. 
Consider lattice points $\{ \bm{t}_i \in [0, 1]^s: i = 0, \dotsc, N-1 \}$ given by~\eqref{eq:lattice_rule} with a generating vector $\bm{z} \in \mathbb{N}^s$ and a random shift $\bm{\Delta} \sim \mathcal{U}([0, 1]^s)$. The RSLR estimator $I_N^{\textrm{lat}} \coloneqq I_N^{\textrm{lat}}(f)$ is then defined as
\begin{align}
	\label{eq:RSLR_estimator}
	I_N^{\textrm{lat}} = \frac{1}{N} \sum_{i=0}^{N-1} f({\bm{t}_i + \bm{\Delta}}).
\end{align}
In the following, we introduce a lemma about the variance of $I_N^{\textrm{lat}}(f)$.
\begin{Lemma}[Variance of RSLR estimator~\cite{lecuyer2000variance}]
	\label{lemma:variance_lattice_rule_rqmc_estimator}
	The variance of the RSLR integration estimator, $\var{I_N^{\mathrm{lat}}}$, is given by
	\begin{equation}
		\label{eq:variance_lattice_rule_rqmc_estimator}
		\var{I_N^{\mathrm{lat}}} = \sum_{\substack{\bm{n} \in \mathbb{Z}^s, \bm{n} \cdot \bm{z} \equiv 0(\textrm{mod } N),\\ \bm{n} \neq 0 }}   \abs*{c_{\bm{n}}}^2,
	\end{equation}
where $c_{\bm{n}}$ are the Fourier coefficients of $f$, given by
\begin{equation*}
	c_{\bm{n}} = \int_{\mathbb{T}^s} f(\bm{t}) e^{-i 2\pi \bm{n} \cdot \bm{t} }d \bm{t},
\end{equation*}
and the set $\{ \bm{n} \in \mathbb{Z}^s: \bm{n} \cdot \bm{z} \equiv 0 \mod N \}$ is known as the dual lattice corresponding to the $N$-point rank-1 lattice rule generated by $\bm{z}$. To summarize: the variance of the RSLR estimator is the summation of the squared magnitudes of Fourier coefficients of the integrand over the dual lattice except the origin. 
\end{Lemma}

\begin{Remark}[Monte Carlo estimator variance]
Let $I_N$ denote the usual Monte Carlo estimator. Then the variance of $I_N$ is given by
\begin{equation}
	\var{I_N} = \frac{1}{N} \sum_{\substack{ \bm{n} \in \mathbb{Z}^s,\\ \bm{n} \neq 0 }}   \abs*{c_{\bm{n}}}^2.
\end{equation}
The variance of the Monte Carlo estimator depends on the sum of the squared magnitudes of the Fourier coefficients of the integrand at all non-zero frequencies. For instance, the reader can refer to~\cite{pilleboue2015variance} for derivations. 
\end{Remark}
\begin{Remark}[Fourier coefficients and integration error]
	The Fourier coefficients play an important role in the integration error of deterministic lattice rules. Specifically,
	consider a lattice rule $\{ \bm{t}_i \in [0, 1]^s: i = 0, \dotsc, N-1 \}$ with a generating vector $\bm{z} \in \mathbb{N}^s$. If $f$ has an absolutely convergent Fourier series, then 
	\begin{equation}
		\label{eq:lattice_rule_l1_integration_error}
		{I - \frac{1}{N}\sum_{i=0}^{N-1} f(\bm{t}_i)} = \sum_{\substack{ \bm{n} \in \mathbb{Z}^s, \bm{n} \cdot \bm{z} \equiv 0(\textrm{mod } N),\\ \bm{n} \neq 0 }} c_{\bm{n}},
	\end{equation}
	where $I = \int_{[0, 1]^s} f(\bm{t})d\bm{t}$. See~\cite{niederreiter1992random,dick2010digital,DKS2013high,dick2022lattice,sloan1994lattice}. 
\end{Remark}
Following~\cite{niederreiter1992random}, we introduce the following definition that characterizes the decay property of Fourier coefficients. 
\begin{Definition}
	\label{definition:niederreiter_space}
	Given $\bm{n} \in \mathbb{Z}^s$, $\alpha>1/2$ and a constant $C < +\infty$, the function class $\mathcal{E}_{\bm{n}}^{\alpha}(C)$ is defined as follows: for $f \in \mathcal{E}_{\bm{n}}^{\alpha}(C)$, its Fourier coefficients $c_{\bm{n}}$ satisfy
	\begin{equation}
		\abs{c_{\bm{n}}} \leq C \prod_{j=1}^{s} {{r_1}{({n}_j)}^{-\alpha}}, \ \bm{n} \neq 0,
	\end{equation}
	where $r_1(n_j) \coloneqq \max{\{1, \abs{n_j}\}}$, for $j = 1, \dotsc, s$. 
\end{Definition}
Notice that in Definition~\ref{definition:niederreiter_space}, we extend the range of $\alpha$ from $\alpha > 1$, as originally assumed in Definition 5.1 of~\cite{niederreiter1992random}, to $\alpha > 1/2$. For $f \in \mathcal{E}_{\bm{n}}^{\alpha}(C)$, the following existence result holds: 
\begin{Lemma}
	There exists a generating vector $\bm{z} \in \mathbb{N}^s$ such that
	\begin{equation}
		\var{I_N^{\mathrm{lat}}} = \mathcal{O}(N^{-2\alpha} (\log N)^{2\alpha s}).
	\end{equation}
\end{Lemma}
Notice that the RSLR variance~\eqref{eq:variance_lattice_rule_rqmc_estimator} and the $L^1$ integration error~\eqref{eq:lattice_rule_l1_integration_error} have similar structures. Thus the proof of the lemma proceeds similarly to the proof of the $L^1$ integration error. For details, refer to~{Section 5.1 of~\cite{DKS2013high}, Chapter 5 of~\cite{niederreiter1992random}, and Chapter 1, Chapter 2 of~\cite{dick2022lattice}}. 

\cite{zaremba1968integration_by_parts} characterized a sufficient condition for a function $f$ to belong to the space $\mathcal{E}_{\bm{n}}^{\alpha}(C)$: Let $\alpha > 1$ be an integer, and assume that for all multi-indices $\bm{k} = (k_1, \dotsc, k_s)$ with $0\leq k_j \leq \alpha - 1, \ j = 1,\dotsc, s$, the mixed partial derivatives
\begin{equation}
	\label{eq:zaremba_condition}
	\frac{\partial^{\abs*{\bm{k}}} f}{\partial t_1^{k_1}\cdots \partial t_s^{k_s}}
\end{equation}
exist and are of bounded variation on $\mathbb{T}^s$ in the sense of Hardy--Krause, where $\abs*{\bm{k}} = k_1 + \cdots + k_s$. Then $f \in \mathcal{E}_{\bm{n}}^{\alpha}(C)$, with the constant $C$ given explicitly. 

In the rest of this section, we establish the relationship between Owen's boundary growth condition and the decay rate $\alpha$ of the Fourier coefficients. 
\begin{Lemma}
	\label{lemma:decay_cn_1d}
	Consider a one-dimensional function $f:[0, 1]\to \mathbb{R}$, whose derivative $\frac{\partial f}{\partial t}$ is defined on $(0, 1)$ and satisfies the boundary growth condition~\eqref{eq:boundary_growth_condition} with $0 \leq A < 1/2$. Then the Fourier coefficients $c_n$ of $f$ satisfy
	\begin{equation}
		\abs{c_n} \lesssim 
		\begin{cases}
			n^{-1+A}, & \text{if}\ 0 < A < 1/2, \\
			n^{-1}\log{n}, & \text{if}\ A = 0.
		\end{cases}
	\end{equation}
\end{Lemma}
The proof, following a classical approach to the Riemann--Lebesgue lemma, is provided to Appendix~\ref{section:proof_decay_cn_1d}. 

When ${A}<0$, we only achieve convergence at the rate of $\mathcal{O}(n^{-1})$ by following the derivations in Appendix~\ref{section:proof_decay_cn_1d}. To achieve higher convergence rates, we propose the following assumptions on higher-order derivatives.
\begin{Assumption}
	\label{assumption:bgc_higher_order_derivative}
	Let $\bar{A} < 1$, $k \in \mathbb{N}_0$. Consider a 1-d function $f: [0, 1]\to \mathbb{R}$. We assume the following conditions:
	\begin{enumerate}
		\item The $(k-1)$-th derivative, $\frac{\partial^{k-1} f}{\partial t^{k-1}}$ is continuous on $\mathbb{T} = [0, 1]$.
		\item The $k$-th derivative $\frac{\partial^{k} f}{\partial t^{k}}$ is piecewise continuous on $\mathbb{T}$ with a finite number of discontinuities. 
		\item The $(k+1)$-th derivative $\frac{\partial^{k+1} f}{\partial t^{k+1}}$ satisfies the boundary growth condition~\eqref{eq:boundary_growth_condition}: $\abs*{\frac{\partial^{k+1} f}{\partial t^{k+1}}} \lesssim \phi(t)^{-1-\bar{A}}$. 
	\end{enumerate}

	When $k = 0$, we define  
	\begin{equation*}
		\frac{\partial^{-1}f}{\partial t^{-1}} (t_0) \coloneqq \int_{0}^{t_0} f(t)dt. 
	\end{equation*}
\end{Assumption}
\begin{Lemma}
	\label{lemma:decay_cn_1d_k_plus_1}
	For a function $f$ satisfying Assumption~\ref{assumption:bgc_higher_order_derivative}, the Fourier coefficients $c_n$ satisfy the following estimates:
		\begin{equation}
			\abs{c_n} \lesssim 
			\begin{cases}
				n^{-1-k}, & \text{if}\  \bar{A} < 0,  \\
				n^{-1-k} \log n, & \text{if}\ \bar{A} = 0, \\
				n^{-1-k+\bar{A}}, & \text{if}\ 0 < \bar{A} < 1.
			\end{cases}
		\end{equation}
\end{Lemma}
\begin{proof}
	We start with the integration by parts.
	\begin{equation}
		\begin{split}
			c_n \cdot (-i2\pi n)^k &= \int_{0}^{1} (-i2\pi n)^k f(t) e^{-i2\pi n t} dt\\
			&= -\int_{0}^{1} (-i2\pi n)^{k-1} \frac{\partial f}{\partial t} e^{-i2\pi n t} dt + \underbrace{(-i 2\pi n)^{k-1} f(t) e^{-i2\pi n t}\mid_{0}^{1}}_{=0}\\
			&\ \vdots\\
			&=(-1)^{k} \int_{0}^{1}  \frac{\partial^{k} f}{\partial t^{k}} e^{-i2\pi n t} dt +  \underbrace{\frac{\partial^{(k-1)} f}{\partial t^{k-1}} e^{-i2\pi n t}\mid_{0}^{1}}_{=0},
		\end{split}
	\end{equation}
	where we have used the continuity of $f$, $\frac{\partial f}{\partial t}$, $\dotsc$, $\frac{\partial^{(k-1)} f}{\partial t^{k-1}}$ on $\mathbb{T}$. Notice that we can use the same shifting technique as in the proof of Lemma~\ref{lemma:decay_cn_1d}. We have
 	\begin{equation}
		\int_{0}^{1}  \frac{\partial^{k} f}{\partial t^{k}}(t) e^{-i2\pi n t} dt = - \int_{0}^{1}  \frac{\partial^{k} f}{\partial t^{k}}\left(t + \frac{1}{2n}\right) e^{-i2\pi n t} dt, 
	\end{equation}
	and
	\begin{equation}
		c_n \cdot (-i2\pi n)^k = \frac{1}{2} \int_{0}^{1} \left( \frac{\partial^{k} f}{\partial t^{k}}(t) - \frac{\partial^{k} f}{\partial t^{k}}\left(t + \frac{1}{2n}\right) \right) e^{-i2\pi n t} dt.
	\end{equation}
	We derive the following upper bound for $\abs{c_n}$:
	\begin{equation}
		\abs{c_n} \lesssim {n^{-k}} \int_{0}^{1} \abs*{ \frac{\partial^{k} f}{\partial t^{k}}(t) - \frac{\partial^{k} f}{\partial t^{k}}\left(t + \frac{1}{2n}\right) }dt. 
	\end{equation}
	{Without loss of generality, assume that the piecewise continuous function $\frac{\partial^{k} f}{\partial t^{k}}$ has finitely many discontinuities at points $x_1, \dotsc, x_d$ where $d < +\infty$ and $0 \leq x_1 < x_2 < \cdots < x_d \leq 1$. } 
	
	We denote the set $\Omega_d = \{{t} \in [0, 1]: \exists i = 1, \dotsc, d, x_i \in [t, t+\frac{1}{2n} ]\}$, i.e. the set of $t \in [0, 1]$ such that the interval $[t, t+\frac{1}{2n} ]$ contains at least one discontinuity point. 
	When $t \in \Omega_d$ and $\bar{A} > 0$, we have
	\begin{equation}
		\label{eq:bound_jump_disc_1}
		\begin{split}
			\int_{\Omega_d} \abs*{ \frac{\partial^{k} f}{\partial t^{k}}(t) - \frac{\partial^{k} f}{\partial t^{k}}\left(t + \frac{1}{2n}\right) } dt &\leq \sum_{i=1}^d \int_{t_i \in [t, t+\frac{1}{2n}]} \abs*{\frac{\partial^{k} f}{\partial t^{k}} (t)} + \abs*{\frac{\partial^{k} f}{\partial t^{k}} (t + \frac{1}{2n})} dt \\
			&\leq 2d \int_{t_i \in [t, t+\frac{1}{2n}]} \abs*{\frac{\partial^{k} f}{\partial t^{k}} (t)} + \abs*{\frac{\partial^{k} f}{\partial t^{k}} (t + \frac{1}{2n})} dt\\
			&\lesssim \int_{[0, \frac{1}{2n}]} \phi(t)^{-\bar{A}} dt.
		\end{split}
	\end{equation}
	Notice that
	\begin{equation}
		\int_{[0, \frac{1}{2n}]} \phi(t)^{-\bar{A}} dt \lesssim \begin{cases}
			n^{-1 + \bar{A}}, & \text{if}\ 0 < \bar{A} < 1, \\
			n^{-1} \log(n), & \text{if}\ \bar{A} = 0.
		\end{cases}
	\end{equation}
	Notice that when $\bar{A} < 0$, the function $\frac{\partial^{k} f}{\partial t^{k}}$ is bounded. Thus, when $t \in \Omega_d$ and $\bar{A} < 0$, we have
	\begin{equation}
		\begin{split}
			\int_{\Omega_d} \abs*{ \frac{\partial^{k} f}{\partial t^{k}}(t) - \frac{\partial^{k} f}{\partial t^{k}}\left(t + \frac{1}{2n}\right) } dt &\lesssim  2d \int_{t_i \in [t, t+\frac{1}{2n}]} \abs*{\frac{\partial^{k} f}{\partial t^{k}} (t)} + \abs*{\frac{\partial^{k} f}{\partial t^{k}} (t + \frac{1}{2n})} dt \lesssim n^{-1}.
		\end{split}
	\end{equation}
	When $t \notin \Omega_d$, following the proof in Appendix~\ref{section:proof_decay_cn_1d}, we have,
	\begin{equation*} 
		\int_{[0, 1] \setminus \Omega_d } \abs*{ \frac{\partial^{k} f}{\partial t^{k}}(t) - \frac{\partial^{k} f}{\partial t^{k}}\left(t + \frac{1}{2n}\right) } dt \leq
		\int_{[0, 1] } \abs*{ \frac{\partial^{k} f}{\partial t^{k}}(t) - \frac{\partial^{k} f}{\partial t^{k}}\left(t + \frac{1}{2n}\right) } dt \lesssim \begin{cases}
			n^{-1 + \bar{A}}, & \text{if}\ 0 < \bar{A} < 1, \\
			n^{-1} \log(n), & \text{if}\ \bar{A} = 0, \\
			n^{-1}, & \text{if}\ \bar{A} < 0.
		\end{cases}
	\end{equation*}
	Thus, we have
	\begin{equation}
		\begin{split}
			\abs{c_n} &\lesssim n^{-k} \left( \int_{[0, 1] \setminus \Omega_d}   \abs*{ \frac{\partial^{k} f}{\partial t^{k}}(t) - \frac{\partial^{k} f}{\partial t^{k}}\left(t + \frac{1}{2n}\right) } dt  + \int_{\Omega_d} \abs*{ \frac{\partial^{k} f}{\partial t^{k}}(t) - \frac{\partial^{k} f}{\partial t^{k}}\left(t + \frac{1}{2n}\right) }  dt \right) \\
			& \lesssim \begin{cases}
				n^{-1 - k + \bar{A}}, & \text{if}\ 0 < \bar{A} < 1, \\
				n^{-1 - k} \log(n), & \text{if}\ \bar{A} = 0, \\
				n^{-1 - k}, & \text{if}\ \bar{A} < 0.
			\end{cases}
		\end{split}
	\end{equation}
	This concludes the proof.
\end{proof}

Under Assumption~\ref{assumption:bgc_higher_order_derivative}, Lemma~\ref{lemma:decay_cn_1d_k_plus_1} extends the result of Lemma~\ref{lemma:decay_cn_1d}. To interpret the result, consider that both the boundary unboundedness and the interior discontinuities can affect the convergence rate. When $\bar{A} \geq 0$, the boundary unboundedness dominates the convergence rate. Conversely, when $\bar{A} < 0$, the $k$-th derivative $\frac{\partial^k f}{\partial t^k}$ is bounded at the boundary, and the convergence rate is dominated by the discontinuities. 
\begin{Remark}[Comparison with Zaremba's condition]
	Compared to Zaremba's sufficient condition~\eqref{eq:zaremba_condition}, which assumes that $\frac{\partial^k f}{\partial t^k}$ is of bounded variation in the Hardy--Krause sense, our approach relaxes this assumption. Specifically, when $\bar{A} = 0$, we assume that $\frac{\partial^k f}{\partial t^k}$ is piecewise continuous but not necessarily of bounded variation in the Hardy--Krause sense, given that $\frac{\partial^k f}{\partial t^k} \lesssim \log (\phi(t))$. Nevertheless, our relaxation leads to a slower convergence rate by a logarithmic factor. 
\end{Remark}

In the following, we study the multidimensional cases, focusing on the first-order mixed partial derivatives of the integrand. The cases involving higher-order mixed partial derivatives can be analyzed similarly using integration by parts and are not detailed in this work. 
\begin{Lemma}
	\label{lemma:decay_cn_multidimension}
	If $f$ satisfies Assumption~\ref{assumption:owen_boundary_growth_condition} with $0 < A_j < 1/2$, $j = 1, \dotsc, s$, then
	\begin{equation}
		\abs{c_{\bm{n}}} \lesssim \prod_{j=1}^{s} r_1(n_j)^{-1+A_j},
	\end{equation}
	with $\bm{n}\in \mathbb{N}_0^s$, $r_1(n_j) \coloneqq \max\{1, \abs{n_j}\}$ for $j = 1, \dotsc, s$. 
\end{Lemma}
\begin{proof}
	Let $\frac{1}{\bm{n}} = (\frac{1}{{n}_1}, \dotsc, \frac{1}{{n}_s})$ for $n_1, \dotsc, n_s \neq 0$. For $\mathfrak{u} \subseteq 1:s$, consider the vector $\bm{v}_{\mathfrak{u}} = \frac{1}{2\bm{n}}^{\mathfrak{u}}:0^{-\mathfrak{u}}$, i.e., $(\bm{v}_{\mathfrak{u}})_j = \frac{1}{2 n_j}$ when $j \in \mathfrak{u}$ and 0 otherwise. We have
	\begin{equation}
		\label{eq:fourier_coefficient_bound_multi_dimension}
		\begin{split}
			c_{\bm{n}} &=\int_{[0, 1]^s} f(\bm{t}) e^{-i2\pi \bm{n} \cdot \bm{t}} d\bm{t} = \int_{[0, 1]^s} f(\bm{t} + \bm{v}_{\mathfrak{u}}) e^{-i2\pi \bm{n} \cdot (\bm{t}+\bm{v}_{\mathfrak{u}})} d\bm{t} = (-1)^{\abs{\mathfrak{u}}} \int_{[0, 1]^s} f(\bm{t} + \bm{v}_{\mathfrak{u}}) e^{-i2\pi \bm{n} \cdot \bm{t}} d\bm{t}.
		\end{split}
	\end{equation}
	By summing over all $2^s$ sets $\mathfrak{u} \subseteq 1:s$ in~\eqref{eq:fourier_coefficient_bound_multi_dimension}, we have
	\begin{equation}
		\begin{split}
			\abs{c_{\bm{n}}} &\leq \frac{1}{2^s} \int_{[0, 1]^s} \abs*{\sum_{\mathfrak{u} \subseteq 1:s} (-1)^{\abs{\mathfrak{u}}} f(\bm{t} + \bm{v}_{\mathfrak{u}})} d\bm{t} := \frac{1}{2^s} \int_{[0, 1]^s} \abs*{\Delta (f;\bm{t}, \bm{t} + \frac{1}{2\bm{n}})} d\bm{t},
		\end{split}
	\end{equation}
	where we use the notation $\Delta (f;\bm{t}, \bm{t} + \frac{1}{2\bm{n}})$ from~\cite{owen2005multidimensional} to denote the alternating sum over the vertices of the hyperrectangle $[\bm{t}, \bm{t} + \frac{1}{2\bm{n}}]$. We proceed the derivation as follows:
	\begin{equation}
		\label{eq:fourier_coefficient_bound_multi_dimension_contd}
		\begin{split}
			\abs*{c_{\bm{n}}} &\leq \frac{1}{2^s}  \int_{[0, 1]^{s}} \abs*{ \Delta(f; \bm{t}, \bm{t} + \frac{1}{2\bm{n}}) } d\bm{t}\\
			&= \frac{1}{2^s}  \int_{[0, 1]^{s}} \abs*{ \int_{[\bm{t}, \bm{t} + \frac{1}{2\bm{n}} \textrm{ mod } 1]} \partial^{1:s} f(\bm{\tau}) d\bm{\tau} } d\bm{t} \\
			&\leq \frac{1}{2^s}  \int_{[0, 1]^{s}} \int_{[\bm{t}, \bm{t} + \frac{1}{2\bm{n}} \textrm{ mod } 1]} \abs*{ \partial^{1:s} f(\bm{\tau}) } d\bm{\tau} d\bm{t}\\
			&\lesssim \frac{1}{2^s}  \int_{[0, 1]^{s}} \int_{[\bm{t}, \bm{t} + \frac{1}{2\bm{n}} \textrm{ mod } 1]} \prod_{j=1}^s \phi({\tau}_j)^{-1-A_j} d\bm{\tau} d\bm{t}\\
			&= \frac{1}{2^s} \prod_{j=1}^s \int_{[0, 1]}   \int_{[{t}_j, {t}_j + \frac{1}{2{n}_j} \textrm{ mod } 1]} \phi({\tau}_j)^{-1-A_j} d{\tau}_j d{t}_j \\
			&\lesssim \prod_{j=1}^s {n}_j^{-1+A_j},
		\end{split}
	\end{equation}
	where notice that the notation $[\bm{t}, \bm{t} + \frac{1}{2\bm{n}} \textrm{ mod } 1]$ from the third line refers to the definition in Section~\ref{sec:notations} and the last line follows by applying similar arguments as in the proof of Appendix~\ref{section:proof_decay_cn_1d}. 
	
	Now we consider the case when some ${n}_j = 0$ for $j = 1, \dotsc, s$. Let $\mathfrak{u}$ be a nonempty proper subset of $1:s$ such that $n_j \neq 0$ if and only if $j \in \mathfrak{u}$. Let $\bm{n} \in \mathbb{N}^s$ and $\bm{n}_{\mathfrak{u}} = \bm{n}^{\mathfrak{u}} : \bm{0}^{-\mathfrak{u}}$. Similarly to the above derivations, we have the following bounds:
	\begin{equation}
		\label{eq:fourier_coefficient_bound_multi_dimension_partial_0}
		\begin{split}
			\abs*{c_{\bm{n}}} &\leq \frac{1}{2^{\abs{\mathfrak{u}}}}  \int_{[0, 1]^{s}} \abs*{ \Delta(f; \bm{t}, \bm{t} + \frac{1}{2\bm{n}_{\mathfrak{u}} }) }d\bm{t}\\
			&= \frac{1}{2^{\abs{\mathfrak{u}}}}  \int_{[0, 1]^{s}} \abs*{ \int_{[\bm{t}_{\mathfrak{u}}, \bm{t}_{\mathfrak{u}} + \frac{1}{2\bm{n}_{\mathfrak{u}}} \textrm{ mod } 1]} \partial^{\mathfrak{u}} f(\bm{\tau}^{\mathfrak{u}} : \bm{t}^{-\mathfrak{u}}) d\bm{\tau}^{\mathfrak{u}} } d\bm{t} \\
			&\leq \frac{1}{2^{\abs{\mathfrak{u}}}}  \int_{[0, 1]^{s}} \int_{[\bm{t}_{\mathfrak{u}}, \bm{t}_{\mathfrak{u}} + \frac{1}{2\bm{n}_{\mathfrak{u}}} \textrm{ mod } 1]} \abs*{ \partial^{\mathfrak{u}} f(\bm{\tau}^{\mathfrak{u}} : \bm{t}^{-\mathfrak{u}}) } d\bm{\tau}^{\mathfrak{u}} d\bm{t}\\
			&\lesssim \int_{[0, 1]^{s}} \int_{[\bm{t}_{\mathfrak{u}}, \bm{t}_{\mathfrak{u}} + \frac{1}{2\bm{n}_{\mathfrak{u}}} \textrm{ mod } 1]} \prod_{j \in \mathfrak{u}} \phi({\tau}_j)^{-1-A_j}  d\bm{\tau}^{\mathfrak{u}} \prod_{k \notin \mathfrak{u}} \phi({t}_k)^{-A_k} d\bm{t}\\
			&= \prod_{j \in \mathfrak{u}} \int_{[0, 1]}  \int_{[{t}_j, {t}_j + \frac{1}{2{n}_j} \textrm{ mod } 1]} \phi({\tau}_j)^{-1-A_j} d{\tau}_j d{t}_j \cdot \prod_{k \notin \mathfrak{u}} \int_{[0, 1]}    \phi({t}_k)^{-A_k}  d{t}_k\\
			&\lesssim \prod_{j\in \mathfrak{u}} {n}_j^{-1+A_j}. 
		\end{split}
	\end{equation}
	By combining~\eqref{eq:fourier_coefficient_bound_multi_dimension_contd} and~\eqref{eq:fourier_coefficient_bound_multi_dimension_partial_0}, we reach the conclusion. 
\end{proof}

The result of Lemma~\ref{lemma:decay_cn_multidimension} can be further extended to functions with axis-parallel discontinuities, as stated in the following lemma. 
\begin{Lemma}
	\label{lemma:decay_cn_multidimension_axis_parallel}
	Let $f: [0, 1]^s \to \mathbb{R}$ be a function of the form $f = \mathbbm{1}_{\Omega} \cdot g$, where $\Omega = [\bm{a}, \bm{b}] \subseteq [0, 1]^s$ is a hyperrectangle. Assume that the first-order mixed partial derivatives $\frac{\partial^{\mathfrak{u}} g}{\partial \bm{t}_{\mathfrak{u}}}$ are defined on $(0, 1)^s$ and satisfy the boundary growth condition~\eqref{eq:boundary_growth_condition} with $0 < A_j < 1/2$ for $j = 1, \dotsc, s$. Then, the Fourier coefficients $c_{\bm{n}}$ satisfy
	\begin{equation}
		\abs{c_{\bm{n}}} \lesssim \prod_{j=1}^{s} r_1(n_j)^{-1+A_j}.
	\end{equation}
\end{Lemma}
The proof can be found in Appendix~\ref{section:proof_decay_cn_multidimension_axis_parallel}. 

In this section, we have established the relationship between the integrand regularity and the Fourier coefficients decay. We found that the Fourier coefficients decay rate is dominated by boundary unboundedness and the interior discontinuities of the function $f$ or its derivatives. Moreover, we demonstrate how Fourier coefficients decay affects the variance of the RSLR estimator. In the next section, we study the Sobol' sequence using spectral analysis. 

\subsection{Sobol' sequence}
Let $N = 2^m$, and define $I_{N}^{\mathrm{sob}} \coloneqq I_{N}^{\mathrm{sob}}(f)$ as the scrambled $(t, m, s)$-Sobol' sequence estimator of the integral of $f$. \cite{owen1997monte, owen1998scrambling} use ANOVA decomposition and the Haar wavelet basis to derive the variance of $I_{N}^{\mathrm{sob}}$. Following~\cite{owen1997monte}, let the multiresolution ANOVA decomposition of $f$ be given by: $f = \sum_{\mathfrak{u} \subseteq 1:s} \sum_{\bm{\ell} \in \mathbb{N}_0^{\abs{\mathfrak{u}}} } f_{\mathfrak{u},\bm{\ell}}$. The variance of $I_{N}^{\mathrm{sob}}$ can be expressed as:
\begin{equation}
	\label{eq:sobol_variance_anova}
	\var{I_{N}^{\mathrm{sob}}} = \frac{1}{N} \sum_{\mathfrak{u} \neq \varnothing} \sum_{\bm{\ell} \in \mathbb{N}_0^{\abs{\mathfrak{u}}}} \Gamma_{\mathfrak{u}, \bm{\ell}} \sigma^2_{\mathfrak{u}, \bm{\ell}},
\end{equation}
where $\sigma^2_{\mathfrak{u}, \bm{\ell}} = \var{f_{\mathfrak{u}, \bm{\ell}}} = \int_{[0, 1]^s} f^2_{\mathfrak{u}, \bm{\ell}}(\bm{t})d\bm{t}$. The term $\Gamma_{\mathfrak{u}, \bm{\ell}}$ is referred to as the ``gain coefficient'' in~\cite{owen1997monte}, which is a quantity depending on the net parameters. Recently, \cite{pan2023nonzero} derived an explicit formula for the gain coefficient in base 2. Additionally, \cite{goda2023improved} derived an upper bound for the gain coefficients in general prime bases. 

\cite{dick2010digital} apply the Walsh basis and derive an equivalent formulation of the variance:
\begin{equation}
	\label{eq:sobol_variance}
	\var{I_N^{\mathrm{sob}}} = \frac{1}{N} \sum_{\bm{\ell} \in \mathbb{N}_0^s} \Gamma_{\bm{\ell}} \sigma^2_{\bm{\ell}}.
\end{equation}
It is worth mentioning that $\Gamma_{\bm{\ell}}$ in~\eqref{eq:sobol_variance} equals $\Gamma_{1:s, \bm{\ell}}$ in~\eqref{eq:sobol_variance_anova}. An alternative formulation of the variance in terms of the sum of squared Walsh coefficients over the dual net is presented in~\cite{l2002recent}. 

In this work, we follow the notations from~\cite{dick2010digital} for simplicity. Next, we study the decay properties of $\sigma^2_{\bm{\ell}}$. We first consider the Walsh decomposition of $f$:
\begin{equation}
	f(\bm{t}) = \sum_{\bm{\ell} \in \mathbb{N}_0^s} \bar{f}({\bm{\ell}}) {}_{2}\wal_{\bm{\ell}}(\bm{t}),
\end{equation}
where $\bar{f}$ refers to the Walsh--Fourier coefficients. The Walsh function ${}_{2}\mathrm{wal}_{\ell}$ is introduced in Appendix~\ref{section:walsh_functions}. Following the notations in~\cite{DKS2013high,dick2010digital}, for $\bm{\ell} = \left({\ell}_1, \dotsc, {\ell}_s \right) \in \mathbb{N}_0^s$, we define the index set
\begin{equation}
	L_{\bm{\ell}} = \{ \bm{k} = ({k}_1, \dotsc, {k}_s) \in \mathbb{N}_0^s: \left\lfloor 2^{{\ell}_j - 1} \right\rfloor \leq k_j < 2^{{\ell}_j}\ \textrm{for}\ 1\leq j \leq s \}.
\end{equation}
Then $\sigma^2_{\bm{\ell}}$ is given by
\begin{equation}
	\label{eq:sigma2_ell}
	\sigma^2_{\bm{\ell}} = \sum_{\bm{k}\in L_{\bm{\ell}}} \abs{\bar{f}(\bm{k})}^2 = \var{\beta_{\bm{\ell}}},
\end{equation}
where the function $\beta_{\bm{\ell}}$ is given by the Walsh expansion of $f$ within the index set $L_{\bm{\ell}}$, i.e.,
\begin{equation}
	\beta_{\bm{\ell}}(\bm{t}) = \sum_{\bm{k}\in L_{\bm{\ell}}} \bar{f}(\bm{k}) {}_{2}\wal_{\bm{k}}(\bm{t}). 
\end{equation}

To proceed with the summations of the Walsh--Fourier coefficients $\bar{f}$ over the set $L_{\bm{\ell}}$, we define another set:
\begin{equation}
	\label{eq:set_T_ell}
	T_{\bm{\ell}} = \{ \bm{k} = ({k}_1, \dotsc, {k}_s) \in \mathbb{N}_0^s: 0 \leq {k}_j < 2^{\ell_j}\ \textrm{for}\ 1\leq j \leq s \}.
\end{equation}
Next, define the operator $\Delta_j$ as set subtraction:
\begin{equation}
	\Delta_j T_{\bm{\ell}} = 
		 T_{\bm{\ell}} \setminus T_{\bm{\ell} - \bm{e}_j},
\end{equation}
where $\bm{e}_j$ is the $j$-th standard unit vector with value 1 in coordinate $j$ and 0 elsewhere. When ${\ell}_j = 0$, we define $T_{\bm{\ell} - \bm{e}_j} \coloneqq \varnothing$. Define $\bigotimes$ as the operator composition, i.e., $\bigotimes_{j=1}^s \Delta_j = \Delta_s \circ \Delta_{s-1} \circ \cdots \circ \Delta_1$. Using this notation, we can express the set $L_{\bm{\ell}}$ as:
\begin{equation}
	L_{\bm{\ell}} = \bigotimes_{j=1}^{s} \Delta_j T_{\bm{\ell}}. 
\end{equation}
We now present a lemma that provides an explicit expression for the Walsh series over the index set $T_{\bm{\ell}}$, which is essential for analyzing the properties of $\sigma^2_{\bm{\ell}}$. 
\begin{Lemma}[Walsh series in $T_{\bm{\ell}}$]
	Given the function $f \in L^2([0, 1]^s)$ and the index set $T_{\bm{\ell}}$ as defined in~\eqref{eq:set_T_ell}, we have the expression for the Walsh series in base 2 over the index set $T_{\bm{\ell}}$ as
	\begin{equation}
		\label{eq:walsh_series_T_ell}
		\sum_{\bm{k} \in T_{\bm{\ell}}} \bar{f}(\bm{k}) {}_{2}\wal_{\bm{k}}(\bm{t}) = \prod_{j=1}^{s} 2^{\ell_j} \int_{\cap_{j=1}^s \left\lfloor  {y}_j 2^{ {\ell}_j}  \right\rfloor = \left\lfloor {t}_j 2^{{\ell}_j}  \right\rfloor} f(\bm{y})d\bm{y}.
	\end{equation}
	\label{lemma:walsh_series_T_ell}
\end{Lemma}
The proof utilizes the Walsh--Dirichlet kernel and is provided in Appendix~\ref{section:walsh_functions}. 

Using Lemma~\ref{lemma:walsh_series_T_ell}, the Walsh series in $L_{\bm{\ell}}$ can be expressed as follows:
	\begin{equation}
		\label{eq:walsh_series_L_ell}
		\begin{split}
			\beta_{\bm{\ell}}(\bm{t}) &= \sum_{\bm{k} \in L_{\bm{\ell}}} \bar{f}(\bm{k}) {}_{2}\wal_{\bm{k}}(\bm{t}) = \sum_{\bm{k} \in \bigotimes_{j=1}^{s} \Delta_j T_{\bm{\ell}}} \bar{f}(\bm{k}) {}_{2}\wal_{\bm{k}}(\bm{t})\\
			&= \bigotimes_{j=1}^{s} \Delta_j \sum_{\bm{k} \in T_{\bm{\ell}}} \bar{f}(\bm{k}) {}_{2}\wal_{\bm{k}}(\bm{t})\\
			&= \bigotimes_{j=1}^{s} \Delta_j \prod_{j=1}^{s} 2^{\ell_j} \int_{\cap_{j=1}^s \left\lfloor {y}_j 2^{{\ell}_j}  \right\rfloor = \left\lfloor {t}_j 2^{{\ell}_j}  \right\rfloor} f(\bm{y})d\bm{y},
		\end{split}
	\end{equation}
	where in the last two equalities we have abused the notation $\Delta_j$ by defining
	\begin{equation}
		\Delta_j \sum_{\bm{k} \in T_{\bm{\ell}}} \bar{f}(\bm{k}) {}_{2}\wal_{\bm{k}}(\bm{t}) =
		\begin{cases} \sum_{\bm{k} \in T_{\bm{\ell}}} \bar{f}(\bm{k}) {}_{2}\wal_{\bm{k}}(\bm{t}) - \sum_{\bm{k} \in T_{\bm{\ell} - \bm{e}_j}} \bar{f}(\bm{k}) {}_{2}\wal_{\bm{k}}(\bm{t}), & \ \mathrm{if} \ {\ell}_j > 1,\\
		\sum_{\bm{k} \in T_{\bm{\ell}}} \bar{f}(\bm{k}) {}_{2}\wal_{\bm{k}}(\bm{t}), & \ \mathrm{otherwise}.
		\end{cases}
	\end{equation}
	We can further express $\beta_{\bm{\ell}}$ in~\eqref{eq:walsh_series_L_ell}. First, we consider the case in 1 dimension. When ${\ell} = 0$, $\beta_{0} = \int_{0}^{1} f(y)dy$. For $\ell \neq 0$, the Equation~\eqref{eq:walsh_series_L_ell} simplifies to
	\begin{equation}
		\label{eq:beta_ell_1d}
		\begin{split}
		\beta_{\ell}(t) &= 2^{\ell} \int_{\left\lfloor y2^{\ell} \right\rfloor = \left\lfloor t2^{\ell} \right\rfloor} f(y) dy - 2^{\ell-1} \int_{\left\lfloor y2^{\ell-1} \right\rfloor = \left\lfloor t2^{\ell-1} \right\rfloor} f(y) dy\\
		&= \begin{cases}
			2^{\ell-1} \left( \int_{\left\lfloor y2^{\ell} \right\rfloor = \left\lfloor t2^{\ell} \right\rfloor} f(y) dy - \int_{\left\lfloor y2^{\ell} \right\rfloor = \left\lfloor t2^{\ell}+1 \right\rfloor} f(y) dy\right), & \text{if}\ \left\lfloor t2^{\ell} \right\rfloor = \left\lfloor t2^{\ell-1} \right\rfloor, \\
			2^{\ell-1} \left( \int_{\left\lfloor y2^{\ell} \right\rfloor = \left\lfloor t2^{\ell} \right\rfloor} f(y) dy - \int_{\left\lfloor y2^{\ell} \right\rfloor = \left\lfloor t2^{\ell}-1 \right\rfloor} f(y) dy\right), & \text{otherwise}
		\end{cases}\\
		&= 2^{\ell-1} \left( \int_{\left\lfloor y2^{\ell} \right\rfloor = \left\lfloor t2^{\ell} \right\rfloor} f(y)  dy - \int_{\left\lfloor y2^{\ell} \right\rfloor = 2\left\lfloor t2^{\ell-1} \right\rfloor + \xi_{\ell}}  f(y)  dy \right),
		\end{split}
	\end{equation}
	where $\xi_{\ell}(t) = \left\lfloor t2^{\ell} \right\rfloor - 2\left\lfloor t2^{\ell-1} \right\rfloor + 1 \textrm{ mod } 2$.
	Notice that the condition $\left\lfloor t2^{\ell} \right\rfloor = \left\lfloor t2^{\ell-1} \right\rfloor$ indicates that the $\ell$-th digit of $t$ in base 2 is 0. Otherwise, the $\ell$-th digit is 1. The last line simplifies the expression by introducing the term $\xi_{\ell}$, which is the complement of $\left\lfloor t2^{\ell} \right\rfloor - 2\left\lfloor t2^{\ell-1} \right\rfloor$ in base 2. 
	
	In the multidimensional case, for $\bm{\ell} \in \mathbb{N}^s$, we derive the following result by extending the one-dimensional case through induction. 
	\begin{equation}
		\label{eq:beta_ell_s}
		\beta_{\bm{\ell}}(\bm{t}) = \prod_{j=1}^s 2^{{\ell}_j - 1} \left( \sum_{\mathfrak{v}\subseteq 1:s} (-1)^{\abs{\mathfrak{v}}} \int_{\cap_{j=1}^s \left\lfloor y_j2^{{\ell}_j} \right\rfloor = 2\left\lfloor t_j2^{{\ell}_j-1} \right\rfloor + \xi_{\bm{\ell}, j} } f(y)dy \right),
	\end{equation}
	with $\xi_{\bm{\ell}, j} = \left\lfloor t2^{{\ell}_j} \right\rfloor - 2\left\lfloor t2^{{\ell}_j-1} \right\rfloor + \mathbbm{1}_{j \in \mathfrak{v}} \textrm{ mod } 2$. When $\bm{\ell} \in \mathbb{N}_0^s$ and some ${\ell}_j \neq 0$ if and only if $j \in \mathfrak{u}$, $\mathfrak{u} \subseteq 1:s$, we have
	\begin{equation}
		\label{eq:beta_ell_partial_s}
		\beta_{\bm{\ell}}(\bm{t}) = \prod_{j \in \mathfrak{u}} 2^{{\ell}_j - 1} \prod_{k \notin \mathfrak{u}} 2^{\bm{\ell}_k} \ \left( \sum_{\mathfrak{v}\subseteq \mathfrak{u}} (-1)^{\abs{\mathfrak{v}}} \int_{[0, 1]^{\abs{-\mathfrak{u}}}} \int_{\cap_{j \in \mathfrak{u}} \left\lfloor \bm{y}_j2^{{\ell}_j} \right\rfloor = 2\left\lfloor \bm{t}_j2^{{\ell}_j-1} \right\rfloor + \xi_{\bm{\ell}, j} } f(\bm{y})d\bm{y}_{\mathfrak{u}}d\bm{y}_{-\mathfrak{u}} \right).
	\end{equation}
	
	We now derive upper bounds for $\sigma^2_{\bm{\ell}}$ using the expressions for $\beta_{\bm{\ell}}$ given in~\eqref{eq:beta_ell_1d}, \eqref{eq:beta_ell_s} and~\eqref{eq:beta_ell_partial_s}. We begin by studying the one-dimensional case in the following lemma. 
	\begin{Lemma}
		For an integrand $f \in L^2([0, 1])$, which satisfies the boundary growth condition~\eqref{eq:boundary_growth_condition} with $-1/2 < A < 1/2$, the term $\sigma^2_{\ell}\coloneqq \sigma^2_{\ell}(f)$, given by~\eqref{eq:sobol_variance}, satisfies the following decay:
		\begin{equation}
			\label{eq:sigma2_ell_decay_1d}
			\sigma^2_{\ell} \lesssim 2^{(2A - 1) \ell}.
		\end{equation}
	\end{Lemma}

\begin{proof}
	Following~\eqref{eq:sigma2_ell} and~\eqref{eq:beta_ell_1d}, for $\ell > 0$, we can express $\sigma^2_{\ell}$ as
	\begin{equation}
		\begin{split}
			\sigma^2_{\ell} &= \int_{[0,1]} \beta^2_{\ell}(x)dx = 2^{2\ell-2} \int_{[0,1]} \left( \int_{\left\lfloor y2^{\ell} \right\rfloor = \left\lfloor t2^{\ell} \right\rfloor} f(y)  dy - \int_{\left\lfloor y2^{\ell} \right\rfloor = 2\left\lfloor t2^{\ell-1} \right\rfloor + \xi_{\ell}}  f(y)  dy \right)^2 dt\\
			&= 2^{2 \ell-2}  \sum_{k = 0}^{2^{\ell}-1} 2^{-\ell} \left( \int_{\left\lfloor y2^{\ell} \right\rfloor = k} f(y) dy - \int_{\left\lfloor y2^{\ell} \right\rfloor =  k+1 } f(y) dy\right)^2 \cdot \mathbbm{1}_{k \textrm{ mod } 2 = 0}\\
			&+ 2^{2 \ell-2}  \sum_{k = 0}^{2^{\ell}-1} 2^{-\ell} \left( \int_{\left\lfloor y2^{\ell} \right\rfloor = k} f(y) dy - \int_{\left\lfloor y2^{\ell} \right\rfloor =  k-1 } f(y) dy\right)^2 \cdot \mathbbm{1}_{k \textrm{ mod } 2 = 1} \\
			&= 2^{2 \ell-2} \sum_{k = 0}^{2^{\ell-1}-1} 2^{-\ell} \cdot 2 \left( \int_{\left\lfloor y2^{\ell} \right\rfloor = 2k} f(y) dy - \int_{\left\lfloor y2^{\ell} \right\rfloor =  2k+1 } f(y) dy\right)^2 \\
			&\leq 2^{\ell-1} \sum_{k = 0}^{2^{\ell-1}-1} \left( \int_{\left\lfloor y2^{\ell} \right\rfloor =  2k } \abs{f(y) - f(y+2^{-\ell})} dy\right)^2. 
		\end{split}
		\label{eq:sigma_2_ell_1d_bound_difference}
	\end{equation}
	We first consider the case $\ell = 1$. When $A\neq 0$, we have
\begin{equation}
	\label{eq:sigma_2_ell_bound_difference_ell_1_A_nonzero}
	\begin{split}
		\sum_{k = 0}^{2^{\ell-1}-1} \left( \int_{\left\lfloor y2^{\ell} \right\rfloor =  2k } \abs{f(y) - f(y+2^{-\ell})} dy \right)^2 &\lesssim \left( \int_{0}^{\frac{1}{2}} \left( \int_{{y}}^{\frac{1}{2}}  {y}_{0}^{-A - 1} d {y}_{0} + \int_{\frac{1}{2}}^{{y} + \frac{1}{2}} (1 - {y}_{0})^{-A - 1} d{y}_{0} \right) d{y} \right)^2 \\
		&= \left( \int_{0}^{\frac{1}{2}} \frac{1}{A} \left[ {y}^{-A} + \left(\frac{1}{2} - {y}\right)^{-A} - 2\left(\frac{1}{2}\right)^{-A} \right] \right)^2\\
		&= \left(\frac{1}{A(1-A)} \right)^2 2^{2A - 2}.
	\end{split}
\end{equation}
When $A = 0$, we have
\begin{equation}
	\label{eq:sigma_2_ell_bound_difference_ell_1_A_zero}
	\begin{split}
		\sum_{k = 0}^{2^{\ell-1}-1} \left( \int_{\left\lfloor y2^{\ell} \right\rfloor =  2k } \abs{f(y) - f(y+2^{-\ell})} dy \right)^2 &= \left( \int_{0}^{\frac{1}{2}} \left( \int_{{y}}^{\frac{1}{2}} {y}_{0}^{ - 1} d{y}_{0} + \int_{\frac{1}{2}}^{{y} + 2^{-1}} (1-{y}_{0})^{- 1} d{y}_{0} \right) d{y} \right)^2 \\
		&= \left( \int_{0}^{\frac{1}{2}}  \left( 2\log\frac{1}{2} - \log ({y}) - \log (\frac{1}{2} - {y}) \right) d {y} \right)^2\\
		&= 1. 
	\end{split}
\end{equation}

	Next, we consider the case $\ell > 1$. First, let $A \neq 0$. For $k = 0,\dotsc, 2^{\ell - 2} - 1$, which correspond to $y_0 \in (0, \frac{1}{2} - 2^{-\ell})$ under the condition $\left\lfloor y2^{\ell} \right\rfloor =  2k$, we have
	\begin{equation}
		\label{eq:sobol_first_half_k}
		\begin{split}
			\abs{f(y_0) - f(y_0+2^{-\ell})} &= \abs*{\int_{y_0}^{y_0+2^{-\ell}} \frac{\partial f}{\partial y} dy} \leq \int_{y_0}^{y_0+2^{-\ell}} \abs*{\frac{\partial f}{\partial y} } dy \lesssim \int_{y_0}^{y_0+2^{-\ell}} y^{-1-A} dy\\
			&= \frac{1}{A} \left[ y_0^{-A} - (y_0+2^{-\ell})^{-A} \right]. 
		\end{split}
	\end{equation}
	For $k = 2^{\ell - 2}, \dotsc, 2^{{\ell} - 1} - 1$, which correspond to the range $y_0 \in (\frac{1}{2}, 1-2^{-\ell})$ under the condition $\left\lfloor y2^{\ell} \right\rfloor =  2k$, we have
	\begin{equation}
		\label{eq:sobol_second_half_k}
		\begin{split}
			\abs{f(y_0) - f(y_0 + 2^{-\ell})} &= \abs*{\int_{y_0}^{y_0 + 2^{-\ell}} \frac{\partial f}{\partial y} dy}
			\leq \int_{y_0}^{y_0 + 2^{-\ell}} \abs*{\frac{\partial f}{\partial y} } dy \lesssim \int_{y_0}^{y_0 + 2^{-\ell}} (1-y)^{-1-A} dy\\
			&= \frac{1}{A} \left[ (1 - y_0 - 2^{-\ell})^{-A} - (1 - y_0)^{-A} \right]. 
		\end{split}
	\end{equation}
	The equations~\eqref{eq:sobol_first_half_k} and~\eqref{eq:sobol_second_half_k} cover all the cases for $k = 0, \dotsc, 2^{\ell - 1} - 1$ in~\eqref{eq:sigma_2_ell_1d_bound_difference}. Given the symmetric results in~\eqref{eq:sobol_first_half_k} and~\eqref{eq:sobol_second_half_k}, it suffices to consider the cases $k = 0, \dotsc, 2^{\ell-2}-1$ in the following derivations. We have
	\begin{equation}
		\begin{split}
			\int_{\left\lfloor y2^{\ell} \right\rfloor =  2k} \abs{f(y) - f(y+2^{-\ell})} dy
			&\lesssim \frac{1}{A} \int_{\left\lfloor y2^{\ell} \right\rfloor =  2k} y^{-A} - (y+2^{-\ell})^{-A} dy\\
			&= \frac{1}{A} \int_{2k\cdot 2^{-\ell}}^{(2k+1)2^{-\ell}} y^{-A} - (y+2^{-\ell})^{-A} dy\\
			&= \frac{1}{A(1-A)} \left[ 2\cdot\left(2^{-\ell}(2k+1)\right)^{1-A} - \left(2^{-\ell}(2k)\right)^{1-A} - \left(2^{-\ell}(2k+2)\right)^{1-A} \right] \\
			&= \frac{2^{-\ell(1-A)}}{A(1-A)} \left[ 2\cdot\left(2k+1\right)^{1-A} - \left(2k\right)^{1-A} - \left(2k+2\right)^{1-A} \right].
		\end{split}
	\end{equation}
	When $k = 0$, $2\cdot\left(2k+1\right)^{1-A} - \left(2k\right)^{1-A} - \left(2k+2\right)^{1-A} = 2 - 2^{1-A}$. For $k \geq 1$, we define the auxiliary function $\varXi:[0, 1] \to \mathbb{R}$, $\varXi(x) = x^{1-A}, A \neq 0$. We apply the Taylor expansion to obtain
	\begin{equation}
		\begin{split}
			(2k)^{1-A} &= (2k+1)^{1-A} - (1-A)(2k+1)^{-A} + \frac{\varXi^{\prime \prime} (2k+1) }{2!} + \cdots + (-1)^m \frac{\varXi^{(m)} (2k+1) }{m!} + \cdots,\\
			(2k+2)^{1-A} &= (2k+1)^{1-A} + (1-A)(2k+1)^{-A} + \frac{\varXi^{\prime \prime} (2k+1) }{2!} + \cdots + (1)^m \frac{\varXi^{(m)} (2k+1) }{m!} + \cdots.
		\end{split}
	\end{equation}
	Then for $k \geq 1$, we have 
	\begin{equation}
		2\cdot\left(2k+1\right)^{1-A} - \left(2k\right)^{1-A} - \left(2k+2\right)^{1-A} = -\sum_{m=1}^{\infty} 2\cdot \frac{\varXi^{(2m)} (2k+1) }{(2m)!}.
	\end{equation}
	Thus, an upper bound for $\sigma^2_{\ell}$ is given by
	\begin{equation}
		\begin{split}
			\sigma^2_{\ell} &\leq 2^{\ell-1} \sum_{k = 0}^{2^{\ell-1}-1} \left( \int_{\left\lfloor y2^{\ell} \right\rfloor =  2k } \abs{f(y) - f(y+2^{-\ell})} dy\right)^2 \\
			& \lesssim 2^{\ell-1} \cdot 2 \sum_{k=0}^{2^{\ell-2} - 1} 2^{-2\ell (1-A)} \left( 2\cdot\left(2k+1\right)^{1-A} - \left(2k\right)^{1-A} - \left(2k+2\right)^{1-A} \right)^2\\
			&\lesssim 2^{\ell(2A-1)} \sum_{k=1}^{2^{\ell-2} - 1} \left(\sum_{m=1}^{\infty}  \frac{-\varXi^{(2m)} (2k+1) }{(2m)!} \right)^2.
		\end{split}
	\end{equation}
	For $A \in (-1/2, 0)\cup(0, 1/2)$, $k \geq 1$, we have
	\begin{equation}
		\begin{split}
			\abs*{\sum_{m=1}^{\infty} \frac{-\varXi^{(2m)} (2k+1)}{(2m)!}} 
			&= \abs*{\sum_{m=1}^{\infty} \frac{(2m-2+A)(2m-1+A)\cdots A(1-A) (2k+1)^{-2m+1-A}}{(2m)!}} \\
			&\lesssim \sum_{m=1}^{\infty} (2k+1)^{-2m+1-A} \\
			&= \frac{(2k+1)^{-1-A}}{1 - (2k+1)^{-2}} \leq {(2k+1)^{-1-A}}.
		\end{split}
	\end{equation}
	Notice that,
	\begin{equation}
		\sum_{k=1}^{2^{\ell - 2} - 1} (2k+1)^{-2-2A} < \sum_{k=1}^{\infty} k^{-2-2A} < +\infty. 
	\end{equation}
	Thus we have,
	\begin{equation}
		\sigma^2_{\ell} \lesssim 2^{(2A - 1)\ell}. 
	\end{equation}
	Now we consider the case $\ell > 1$ and $A = 0$. From~\eqref{eq:sigma_2_ell_1d_bound_difference}, we have
	\begin{equation}
		\begin{split}
			\sigma^2_{\ell} &\leq 2^{\ell-1} \sum_{k=0}^{2^{\ell - 1} - 1} \left( \int_{\left\lfloor y2^{\ell} \right\rfloor =  2k } \abs{f(y) - f(y+2^{-\ell})} dy\right)^2\\
			&\lesssim 2^{\ell} \sum_{k=0}^{2^{\ell - 2} - 1} \left( \int_{2k\cdot 2^{-\ell}}^{(2k+1)2^{-\ell}} \log(y+2^{-\ell}) - \log(y) dy \right)^2\\
			&= 2^{\ell} \sum_{k=0}^{2^{\ell - 2} - 1} \left( \varTheta\left((2k+2)2^{-\ell}\right) + \varTheta(2k \cdot 2^{-\ell}) - 2\varTheta((2k+1) \cdot 2^{-\ell}) \right)^2,
		\end{split}
	\end{equation}
	where in the second line we use the symmetry of the boundary growth condition and in the third line we define the auxiliary function $\varTheta:[0, 1]\to \mathbb{R}$, $\varTheta(x) = x\log x - x$. The derivatives of $\varTheta$ are given by $\varTheta^{(2m)}(2k+1) = (2m-2)! (2k+1)^{-2m+1}$, for $m = 1, \dotsc, \infty$. For $k \in \mathbb{N}_0$, similar to earlier derivations we have
	\begin{equation}
		\begin{split}
			\varTheta\left((2k+2)2^{-\ell}\right) &= \varTheta\left((2k+1)2^{-\ell}\right) + 2^{-\ell} \varTheta^{\prime}\left((2k+1)2^{-\ell}\right) + \cdots + (2^{-\ell})^m \frac{\varTheta^{(m)} ((2k+1)2^{-\ell}) }{m!} + \cdots\\
			\varTheta\left((2k)2^{-\ell}\right) &= \varTheta\left((2k+1)2^{-\ell}\right) - 2^{-\ell} \varTheta^{\prime}\left((2k+1)2^{-\ell}\right) + \cdots + (-2^{-\ell})^m \frac{\varTheta^{(m)} ((2k+1)2^{-\ell}) }{m!} + \cdots,
		\end{split}
	\end{equation}
	and
	\begin{equation}
		\begin{split}
			\varTheta\left((2k+2)2^{-\ell}\right) + \varTheta(2k \cdot 2^{-\ell}) - 2\varTheta((2k+1) \cdot 2^{-\ell}) &= \sum_{m=1}^{\infty} 2 (2^{-\ell})^{2m} \frac{\varTheta^{(2m)} ((2k+1)2^{-\ell}) }{(2m)!}\\
			&= 2 \sum_{m=1}^{\infty} (2^{-\ell})^{2m} \frac{(2m-2)! ((2k+1)2^{-\ell})^{-2m+1} }{(2m)!}\\
			&= 2^{-\ell + 1} \sum_{m=1}^{\infty} \frac{(2m-2)! (2k+1)^{-2m+1} }{(2m)!}\\
			& \lesssim 2^{-\ell + 1} (2k+1)^{-1}.
		\end{split}
	\end{equation}
	Finally, for $\ell > 1$ and $A = 0$, we have
	\begin{equation}
		\label{eq:sigma_2_ell_bound_difference_A_0}
			\sigma^2_{\ell} \lesssim 2^{\ell} \sum_{k=0}^{2^{\ell - 2} - 1} 2^{-2\ell} (2k+1)^{-2} \leq 2^{-\ell} \sum_{k=0}^{\infty} (2k+1)^{-2} \lesssim 2^{-\ell}.
	\end{equation}
	This concludes the proof. 
\end{proof}

We present the following lemma for the multidimensional case.
\begin{Lemma}
	\label{lemma:variance_sobol_multi_d}
	We consider an integrand $f \in L^2([0, 1]^s)$ that satisfies the boundary growth condition~\eqref{eq:boundary_growth_condition} with $-1/2 < A^* < 1/2$, where $A^* = \max_{j=1}^s A_j$. Given $\bm{\ell} \in \mathbb{N}_0^s$ with ${\ell}_j \neq 0$ if and only if $j \in \mathfrak{u} \subsetneq 1:s$, the term $\sigma^2_{\bm{\ell}}\coloneqq \sigma^2_{\bm{\ell}}(f)$, given by~\eqref{eq:sobol_variance}, satisfies the following decay:
	\begin{equation}
		\sigma^2_{\bm{\ell}} \lesssim 2^{(2A^* - 1) \abs{\bm{\ell}} }.
	\end{equation}
\end{Lemma}
\begin{proof}
Following~\eqref{eq:beta_ell_partial_s}, we have
\begin{equation}
	\label{eq:sigma2_ell_multi_dimension}
	\begin{split}
		\sigma^2_{\bm{\ell}} &= \int_{[0,1]^s} \abs*{\beta_{\bm{\ell}}(\bm{t})}^2 d\bm{t}\\
		&= \sum_{\substack{\bm{k}_{\mathfrak{u}} \in \mathbb{N}_0^{\abs{\mathfrak{u}}} \\ {\bm{k}_{\mathfrak{u}}\leq 2^{\bar{\bm{\ell}}_{\mathfrak{u}} } - 1} }} \left(\int_{[0, 1]^{\abs{-\mathfrak{u}}} } \int_{\cap_{j \in \mathfrak{u}} \left\lfloor {y}_j 2^{{\ell}_j} \right\rfloor = 2\bm{k}_j}  \sum_{\mathfrak{v} \subseteq \mathfrak{u}} (-1)^{\abs{\mathfrak{v}}} f\left(\bm{y}^{\mathfrak{v}}:\left(\bm{y} + 2^{-\bm{\ell}} \right)^{\mathfrak{u}-\mathfrak{v}} \mid \bm{y}_{-\mathfrak{u}} \right) d\bm{y}_{\mathfrak{u}} d\bm{y}_{-\mathfrak{u}} \right)^2\\
		&= \sum_{\substack{\bm{k}_{\mathfrak{u}} \in \mathbb{N}_0^{\abs{\mathfrak{u}}} \\ {\bm{k}_{\mathfrak{u}}\leq 2^{\bar{\bm{\ell}}_{\mathfrak{u}} } - 1} }} \left( \int_{(E_{\bm{\ell}, 2\bm{k}})_{\mathfrak{u}}} 
		\sum_{\mathfrak{v} \subseteq \mathfrak{u}} (-1)^{\abs{\mathfrak{v}}} \int_{[0, 1]^{\abs{-\mathfrak{u}}} } f\left(\bm{y}^{\mathfrak{v}}:\left(\bm{y} + 2^{-\bm{\ell}} \right)^{\mathfrak{u}-\mathfrak{v}} \mid \bm{y}_{-\mathfrak{u}} \right) d\bm{y}_{-\mathfrak{u}} d\bm{y}_{\mathfrak{u}}  \right)^2\\
		&= \sum_{\substack{\bm{k}_{\mathfrak{u}} \in \mathbb{N}_0^{\abs{\mathfrak{u}}} \\ {\bm{k}_{\mathfrak{u}}\leq 2^{\bar{\bm{\ell}}_{\mathfrak{u}} } - 1} }} \left( \int_{{(E_{\bm{\ell}, 2\bm{k}})_{\mathfrak{u}}}} \int_{[\bm{y}, \bm{y} + 2^{-\bm{\ell}_{\mathfrak{u}}}] } \partial^{\mathfrak{u}} \int_{[0, 1]^{\abs{-\mathfrak{u}}}} f(\bm{y}_0 \mid \bm{y}_{-\mathfrak{u}}) d\bm{y}_{-\mathfrak{u}} d\bm{y}_0 d\bm{y}_{\mathfrak{u}} \right)^2\\ 
		&\lesssim  \sum_{\substack{\bm{k}_{\mathfrak{u}} \in \mathbb{N}_0^{\abs{\mathfrak{u}}} \\ {\bm{k}_{\mathfrak{u}}\leq 2^{\bar{\bm{\ell}}_{\mathfrak{u}} } - 1} }} \left(  \int_{{(E_{\bm{\ell}, 2\bm{k}})_{\mathfrak{u}}}}
		\int_{[\bm{y}, \bm{y} + 2^{-\bm{\ell}_{\mathfrak{u}}}] } \prod_{j\in \mathfrak{u}} \phi({y}_{0, j})^{-A_j - 1} d\bm{y}_0 d\bm{y}_{\mathfrak{u}} \right)^2, 
	\end{split}
\end{equation}
where we denote $\bar{\bm{\ell}}_{\mathfrak{u}} = \bm{\ell}_{\mathfrak{u}} - \bm{1}$, with the subtraction taken component-wise. 
We have
\begin{equation}
	\begin{split}
		& \sum_{\substack{\bm{k}_{\mathfrak{u}} \in \mathbb{N}_0^{\abs{\mathfrak{u}}} \\ {\bm{k}_{\mathfrak{u}}\leq 2^{\bar{\bm{\ell}}_{\mathfrak{u}} } - 1} }} \left(  \int_{{(E_{\bm{\ell}, 2\bm{k}})_{\mathfrak{u}}}} \int_{[\bm{y}, \bm{y} + 2^{-\bm{\ell}_{\mathfrak{u}}}] } \prod_{j\in \mathfrak{u}} \phi({y}_{0, j})^{-A_j - 1} d\bm{y}_0 d\bm{y}_{\mathfrak{u}} \right)^2 \\
		&= \prod_{j\in \mathfrak{u}} \sum_{\substack{{k}_{j} \in \mathbb{N}_0 \\ {{k}_{j}\leq 2^{{{\ell}}_{j} - 1} - 1} }}  \left(  \int_{ \lfloor y_j 2^{\ell_j} \rfloor = 2k_j }
		\int_{[{y}_j, {y}_j + 2^{-{\ell}_{j}}] }  \phi({y}_{0, j})^{-A_j - 1} d {y}_{0,j} d {y}_{j} \right)^2.
	\end{split}
\end{equation}
Refer to Equations~\eqref{eq:sigma_2_ell_1d_bound_difference}-\eqref{eq:sigma_2_ell_bound_difference_A_0}, we have
\begin{equation}
	\label{eq:sigma_2_ell_bound_2_2A_1}
	\sigma_{\bm{\ell}}^2 \lesssim \prod_{j=1}^s 2^{(2A_j - 1){\ell}_j}  \lesssim 2^{(2A^{*} - 1)\abs{\bm{\ell}}}. 
\end{equation}
This concludes the proof. 
\end{proof}
\begin{Theorem}
	Let $f \in L^2[0, 1]^s$ satisfy the boundary growth condition~\eqref{eq:boundary_growth_condition} with $-1/2 < A^{*} < 1/2$, the variance of the RQMC estimator with scrambled Sobol' sequence $I_N^{\mathrm{sob}}$ satisfies
	\begin{equation}
		\var{I_N^{\mathrm{sob}}} = \mathcal{O}(N^{2A^{*} - 2} \log(N)^{s-1}).
	\end{equation}
\end{Theorem}
\begin{proof}
	From Lemma~\ref{lemma:variance_sobol_multi_d}, we obtain $\sigma^2_{\bm{\ell}} \lesssim 2^{(2A^{*} - 1)\abs{\bm{\ell}}}$. The rest of the proof follows the arguments Theorem 3 of~\cite{owen2008local} and Theorem 13.25 of~\cite{dick2010digital}.
\end{proof}

Notice that in our proof, the convergence rate of the Sobol' sequence does not exceed $\mathcal{O}(N^{-3/2})$, even if the integrand has better regularity. This observation aligns with the remarks of~\cite{dick2010digital}, who introduce higher-order digital nets to exploit faster convergence rate. In contrast, the convergence rate of the RSLR largely depends on the integrand regularity.

Moreover, for an integrand satisfying the boundary growth condition~\eqref{eq:boundary_growth_condition} with $-1/2 <  A^{*} < -1/2$, the optimal convergence rate that we obtain is $\mathcal{O}(N^{-3+\delta} \log(N)^{s-1})$, for any $\delta > 0$. This rate closely approximates $\mathcal{O}(N^{-3} (\log N)^{s-1})$, as established in~\cite{yue1999variance, owen2008local} for Lipschitz continuous integrands. However, our analysis does not require the Lipschitz condition.

\section{Unbounded integrands with interior discontinuities}
\label{sec:discontinuous_integrand}
In this section, we aim to extend the analysis of the Sobol' sequence to discontinuous integrands along with the unboundedness as characterized in the boundary growth condition~\eqref{eq:boundary_growth_condition} outlined in Section~\ref{sec:background}. Specifically, following~\cite{he2015convergence}, we consider an integrand of the form $f(\bm{y}) = g(\bm{y}) \cdot \mathbbm{1}_{\bm{y} \in \Omega}$, where $\Omega \subset [0, 1]^s$. Moreover, we assume that $\Omega$ has a partially axis-parallel structure, i.e.,
\begin{equation}
	\label{eq:partially_axis_parallel_sets}
	\Omega = \Omega_{\mathfrak{u}} \times \prod_{j \notin \mathfrak{u}} [a_j, b_j],
\end{equation}
where $\mathfrak{u} \subseteq 1:s$, $\Omega_{\mathfrak{u}}$ is a Lebesgue measurable subset of $\prod_{j \in \mathfrak{u}}[0, 1]$, and $0\leq a_j < b_j\leq 1$ for all $j \notin \mathfrak{u}$. Following~\cite{he2015convergence}, we introduce the following condition on the set $\Omega$. 
\begin{Definition}[Minkowski content~\cite{ambrosio2008outer}]
	Let $(\partial \Omega)_{\epsilon}$ be the outer parallel body of the boundary $\partial \Omega$ of $\Omega$ at a distance $\epsilon$. If the limit
	\begin{equation*}
		\mathcal{M}(\partial \Omega) = \lim_{\epsilon \to 0} \frac{\lambda_s((\partial \Omega)_{\epsilon})}{\epsilon},
	\end{equation*}
	exists and is finite, then $\Omega$ is said to admit an $(s-1)$-dimensional Minkowski content with respect to the s-dimensional Lebesgue measure $\lambda_s$. 
\end{Definition}

As discussed in~\cite{he2015convergence}, when $\Omega_{\mathfrak{u}}$ admits a $\abs{\mathfrak{u}} - 1$-dimensional Minkowski content, the convergence rate of the RQMC estimator for integrands of the form $f = \mathbbm{1}_{\Omega} g$ is given by $\mathcal{O}(n^{-\frac{1}{2} - \frac{1}{4 \abs{\mathfrak{u}} - 2}} (\log n)^{\frac{2s}{2\abs{\mathfrak{u}} - 1}} )$, if $g$ has bounded Hardy--Krause variation. The result improves to $\mathcal{O}(n^{-\frac{1}{2} - \frac{1}{2 \abs{\mathfrak{u}}}})$ if $g \equiv 1$. In a further exploration~\cite{he2018quasi_discontinuous_boundary}, the author considers the case when $g$ is unbounded and satisfies Owen's boundary growth condition~\eqref{eq:boundary_growth_condition}, which is exactly the scenario we consider in this section. The author shows that the convergence rate of the RQMC estimator with the Sobol' sequence is $\mathcal{O}(n^{(A^*-1)(\frac{1}{2} + \frac{1}{4 \abs{\mathfrak{u}} - 2})} (\log n)^{\frac{(1-A^*)s}{2\abs{\mathfrak{u}} - 1}} )$. This results suggests a multiplicative effect arising from both the discontinuity and the boundary unboundedness. In the subsequent sections, we perform a spectral analysis for such discontinuous integrands. 

We start by presenting two lemmas that will be used in the following sections. The first lemma is an excerpt from Lemma 3.4 in~\cite{he2015convergence}. 
\begin{Lemma}[Excerpt from Lemma 3.4 in~\cite{he2015convergence}]
	\label{lemma:elementary_intervals_boundary}
	Consider a vector $\bm{r} \in \mathbb{N}_0^s$ with identical components, where ${r}_j = r$ for all $j = 1, \dotsc, s$, $r \in \mathbb{N}$. We define the elementary interval $E_{\bm{r}, \bm{k}}$ as:
	\begin{equation}
		\label{eq:elementary_interval_r}
		E_{\bm{r}, \bm{k}} = \prod_{j=1}^s \left[ \frac{ {k}_j}{{r}}, \frac{{k}_j + 1}{r} \right), \quad 0 \leq {k}_j \leq r - 1.
	\end{equation}
	Let ${\mathcal{T}_{\mathrm{tot}}^{r} }$ denote the set of the elementary intervals $E_{\bm{r}, \bm{k}}$ that intersect with $\Omega$, and let ${\mathcal{T}_{\mathrm{bdy}}^{r} }$ denote the set of elementary intervals $E_{\bm{r}, \bm{k}}$ that intersect with $\partial \Omega$. The cardinality of ${\mathcal{T}_{\mathrm{bdy}}^{r}}$ has the following upper bound:
	\begin{equation}
		\label{eq:elementary_intervals_boundary_count_bound}
		\abs*{\mathcal{T}_{\mathrm{bdy}}^{r}} \lesssim r^{s-1}.
	\end{equation}
	Furthermore, ${\mathcal{T}_{\mathrm{tot}}^{r}}$ can be expressed as the union of at most $\abs*{\mathcal{T}_{\mathrm{bdy}}^{r}} \lesssim r^{s-1}$ axis-parallel sets of the form $[\bm{a}, \bm{b})$. 
\end{Lemma}
Notice that in Lemma~\ref{lemma:elementary_intervals_boundary}, the elementary intervals $E_{\bm{r}, \bm{k}}$ are isotropic, with a length of $1/r$ in each dimension. This differs from the elementary intervals $E_{\bm{\ell}, \bm{k}}$ in~\eqref{eq:elementary_interval}. 
 
The second Lemma discusses integration within each of the elementary interval. 
\begin{Lemma}
	\label{lemma:alternating_sum_integration}
	Let $f = \mathbbm{1}_{\Omega} \cdot g$, where $g$ satisfies the boundary growth condition~\eqref{eq:boundary_growth_condition} with $A^{*} < \frac{1}{2}$ and $\Omega \subset [0, 1]^s$. For ${\bm{\ell}} \in \mathbb{N}_0^s$, define $\bar{\bm{\ell}} \coloneqq \max(\bm{\ell} - \bm{1}, 0)$, where $\bm{1}$ is the vector of 1's and both subtraction and maximum are taken component-wise. We also define $f(\bm{y}) \coloneqq 0$ when $\bm{y} \notin [0, 1)^s$ for the notation simplicity in the derivations. Then, for all $\bm{k} \in \mathbb{N}_0^s$ such that $0 \leq {k}_j \leq 2^{\bar{{\ell}}_j} - 1$ for $j = 1, \dotsc, s$, the integration of the alternating sum $\Delta(f;\bm{y},\bm{y}+2^{-\bm{\ell}})$ over the elementary interval $E_{{\bm{\ell}}, 2 \bm{k}}$ have the following upper bound:
	\begin{equation}
		\left(\int_{E_{\bm{\ell},2 \bm{k}}} \Delta(f;\bm{y},\bm{y}+2^{-\bm{\ell}}) d\bm{y}\right)^2 \lesssim 
			2^{-\abs{\bm{\ell}}}.
	\end{equation}
\end{Lemma}

\begin{proof}
Notice that we have the following inequality:
\begin{equation}
		\left(\int_{E_{\bm{\ell}, 2\bm{k}}} \Delta(f;\bm{y},\bm{y}+2^{-\bm{\ell}}) d\bm{y}\right)^2 \leq \left(\int_{E_{\bm{\ell}, 2 \bm{k}}} \abs*{\Delta(f;\bm{y},\bm{y}+2^{-\bm{\ell}})} d\bm{y}\right)^2 \leq \left(\int_{E_{\bar{\bm{\ell}}, \bm{k}}} \abs*{f(\bm{y})} d\bm{y}\right)^2.
\end{equation}
When $A^{*} > 0$, the integration in each elementary interval $E_{\bar{\bm{\ell}}, \bm{k}}$ has the following upper bound:
\begin{equation}
\begin{split}
	\left(\int_{E_{\bar{\bm{\ell}}, \bm{k}}} \abs*{f(\bm{y})} d\bm{y}\right)^2 &\lesssim \left( \int_{[0, 2^{-\bar{\bm{\ell}} }]} \bm{y}^{-A^*} d\bm{y} \right)^2  = \left( \prod_{j=1}^s 2^{-\bar{\bm{\ell}}_{j} (1 - A^{*}) } \right)^2 \leq 2^{-\abs{\bar{\bm{\ell}}}} \leq 2^{-s} 2^{-\abs{ {\bm{\ell}}}}.
\end{split}	
\end{equation}
When $A^{*} = 0$, we have
\begin{equation}
\begin{split}
	\left(\int_{E_{\bar{\bm{\ell}}, \bm{k}}} \abs*{f(\bm{y})} d\bm{y}\right)^2 &\lesssim \left( \int_{[0, 2^{-\bar{\bm{\ell}} }]} -\log (\bm{y})  d\bm{y} \right)^2 \lesssim 2^{-2\abs{\bar{\bm{\ell}}}} \prod_{j=1}^s \left( 1 + \bar{\bm{\ell}}_j \right)^2 \lesssim 2^{-\abs{\bm{\ell}}}.
\end{split}
\end{equation}
Otherwise $\abs*{f(\bm{y})}$ is bounded in each elementary interval and
\begin{equation}
\begin{split}
	\left( \int_{E_{\bar{\bm{\ell}}, \bm{k}}} \abs*{f(\bm{y})} d\bm{y} \right)^2 &\leq  \int_{E_{\bar{\bm{\ell}}, \bm{k}}}  f^2(\bm{y})  d\bm{y} \leq \int_{[0, 2^{-\bar{\bm{\ell}} }]} 1 d\bm{y}  \leq 2^{-s} 2^{-\abs{\bm{\ell}}}.
\end{split}	
\end{equation}
This completes the proof.
\end{proof}

With these ingredients, we proceed to analyze various cases for $\Omega$ involved in the discontinuous integrands in the following sections.  

\subsection{Axis-parallel sets}
We first examine the case of axis-parallel sets. 
\begin{Lemma}
	\label{lemma:sobol_axis_parallel}
	Consider a function $g: [0, 1]^s \to \mathbb{R}$ satisfying the boundary growth condition~\eqref{eq:boundary_growth_condition} with $-1/2 < A^{*} < 1/2$. Let $f = \mathbbm{1}_{\Omega} \cdot g$, where $\Omega$ is given in the axis-parallel form $[\bm{a}, \bm{b}) \subset [0, 1]^s$. Then the term $\sigma^2_{\bm{\ell}}\coloneqq \sigma^2_{\bm{\ell}}(f)$, $\bm{\ell} \in \mathbb{N}_0^s$, as defined by~\eqref{eq:sobol_variance}, satisfies the following decay:
	\begin{equation}
		\label{eq:convergence_sobol_axis_parallel}
		\sigma^2_{\bm{\ell}} \lesssim 2^{\max(-1+2A^{*}, -1){\abs{\bm{\ell}}}}. 
	\end{equation}
\end{Lemma}
The proof is contained as a special case in the proof of the next lemma in Appendix~\ref{sec:proof_sobol_axis_parallel_multiple}. 

Next, we introduce a lemma concerning the case where $\Omega$ is the union of multiple axis-parallel sets, which will be utilized in the subsequent proofs. 
\begin{Lemma}
	\label{lemma:sobol_axis_parallel_multiple}
	Consider a function $g: [0, 1]^s \to \mathbb{R}$ that satisfies the boundary growth condition~\eqref{eq:boundary_growth_condition} with $-1/2 < A^{*} < 1/2$. Let $f_M = \mathbbm{1}_{ {\Omega}_M} \cdot g$, where $\Omega_M$ is the union of $M$ intervals of the form $[\bm{a}, \bm{b})$. Then the variance of the scrambled Sobol' sequence estimator satisfies:
	\begin{equation}
		\label{eq:convergence_sobol_axis_parallel_multi_no_sign_change}
		\var{I_N^{\mathrm{sob} }(f_M)} \lesssim M \left(  N^{-2+2A^{*}} \log(N)^{s-1} + N^{-2 } \log(N)^{s-1} \right). 
	\end{equation}
\end{Lemma}
The proof is deferred to Appendix~\ref{sec:proof_sobol_axis_parallel_multiple}. 

\subsection{General discontinuous integrand}
\label{sec:general_discontinuous_integrand_sobol}
We consider a general discontinuous integrand, $f = \mathbbm{1}_{\Omega} \cdot g$, where $\Omega$ admits an $s-1$ dimensional Minkowski content. 

For $\bm{r}$ and $r$ defined in~\eqref{eq:elementary_interval_r}, recall the definition ${\mathcal{T}_{\mathrm{tot}}^{r} } \coloneqq \{ \bigcup E_{\bm{r}, \bm{k}} : E_{\bm{r}, \bm{k}} \cap \Omega \neq \varnothing, 0 \leq {k}_j \leq {r} - 1 , j \in 1:s\}$. 
Define $f_{r} \coloneqq \mathbbm{1}_{{\mathcal{T}_{\mathrm{tot}}^{r} }} g$. Then we have the following inequality:
\begin{equation}
	\var{I_N^{\mathrm{sob}} (f)} \leq 2 \var{I_N^{\mathrm{sob}}(f - f_{r})} + 2 \var{I_N^{\mathrm{sob}}(f_{r})}. 
\end{equation}

From Lemma~\ref{lemma:elementary_intervals_boundary}, ${\mathcal{T}_{\mathrm{tot}}^{r} }$ can be expressed as the union of at most $r^{s-1}$ axis-parallel sets. From Lemma~\ref{lemma:sobol_axis_parallel_multiple}, for a fixed $r \in \mathbb{N}$, we have that 
\begin{equation}
	\var{I_N^{\mathrm{sob}}(f_{r})} \lesssim r^{s-1}  \left(N^{-2} (\log N)^{s-1} +  N^{2A^{*} - 2} (\log N)^{s-1}\right).
\end{equation}
Notice that the support of $f - f_r$ is contained in ${\mathcal{T}_{\mathrm{bdy}}^{r} }$. For a bounded $f$, i.e., when $A^{*} < 0$, we have
\begin{equation}
	\E{(f - f_{r})^2} \lesssim r^{s-1} r^{-s} = r^{-1}. 
\end{equation}
Then we obtain
\begin{equation}
	\begin{split}
		\var{I_N^{\mathrm{sob}}(f - f_{r})} &\leq \frac{1}{N} \var{f - f_{r}} \leq \frac{1}{N} \E{(f - f_{r})^2} \lesssim r^{-1}N^{-1},
	\end{split}
\end{equation}
and consequently,
\begin{equation}
	\begin{split}
		\var{I_N^{\mathrm{sob}} (f)} &\lesssim r^{-1} N^{-1} + r^{s-1} N^{-2} (\log N)^{s-1} + r^{s-1} N^{2A^{*} - 2} (\log N)^{s-1}\\
		&\lesssim r^{-1} N^{-1} + r^{s-1} N^{-2} (\log N)^{s-1}.
	\end{split}
\end{equation}
We choose $r = N^{\frac{1}{s}} (\log N)^{\frac{1}{s} - 1}$ to balance the two terms, yielding
\begin{equation}
	\var{I_N^{\mathrm{sob}} (f)} \lesssim N^{-1-\frac{1}{s}} \left( \log N \right)^{1 - \frac{1}{s}}.
\end{equation}	
When $A^{*} \geq 0$, we apply the H\"older inequality to obtain
\begin{equation}
	\label{eq:holder_inequality}
	\begin{split}
		\E{(f - f_{r})^2} &= \int_{[0, 1]^s } f(\bm{y})^2 \mathbbm{1}_{{\mathcal{T}_{\mathrm{tot}}^{r} } \setminus \Omega} d\bm{y} \leq \int_{[0, 1]^s } f(\bm{y})^2  \mathbbm{1}_{\mathcal{T}_{\mathrm{bdy}}^{r} } d\bm{y}\\
		& \leq \left( \int_{[0, 1]^s } f(\bm{y})^{2p} \right)^{1/p} \left( \int_{[0, 1]^s }   \mathbbm{1}_{\mathcal{T}_{\mathrm{bdy}}^{r} }^{q} d\bm{y}\right)^{1/q},
	\end{split}
\end{equation}
where $1/p + 1/q = 1$. We choose $2A^{*} p < 1$ to ensure the integrability of $f^{2p}$, and consequently $1/q = 1 - 2A^{*} - \delta$ for $0 < \delta < 1 - 2A^{*}$. We have
\begin{equation}
	\begin{split}
		\var{I_N^{\mathrm{sob}}(f - f_{r})} &\leq \frac{1}{N} \E{(f - f_{r})^2} \leq \frac{1}{N} \left(  \int_{[0, 1]^s }   \mathbbm{1}_{\mathcal{T}_{\mathrm{bdy}}^{r} }^{q} d\bm{y} d\bm{y}\right)^{1 - 2A^{*} - \delta}\\
		&\lesssim \frac{1}{N} (r^{s-1} r^{-s})^{1 - 2A^{*} - \delta}\\
		& = r^{-1 + 2A^{*} + \delta} N^{-1},
	\end{split}
\end{equation}
and
\begin{equation}
	\var{I_N^{\mathrm{sob}} (f)} \lesssim r^{-1 + 2A^{*} + \delta} N^{-1} + r^{s-1}  N^{2A^{*} - 2} (\log N)^{s-1}.
\end{equation}
We choose $r = N^{\frac{1 - 2A^{*} }{s - 2A^{*} -\delta }} (\log N)^{\frac{1-s}{s - 2A^{*} - \delta} }$ to balance the two terms and obtain
\begin{equation}
	\begin{split}
		\var{I_N^{\mathrm{sob}} (f)} &\lesssim N^{-1 + \frac{(1-2A^{*})(-1+2A^{*}+\delta)}{s-2A^{*}-\delta}} \left( \log N \right)^{\frac{(s-1)(1-2A^{*} -\delta)}{s-2A^{*} -\delta }}. 
	\end{split}
\end{equation}
The results are summarized in the following Lemma. 
\begin{Lemma}
	For $f = g \cdot \mathbbm{1}_{\Omega}$, where $g$ satisfies the boundary growth condition~\eqref{eq:boundary_growth_condition} with $-1/2 < A^{*} < 1/2$ and $\Omega$ admits an $s-1$ dimensional Minkowski content, the variance of the RQMC estimator with scrambled Sobol' sequence, $I_N^{\mathrm{sob}} \coloneqq I_N^{\mathrm{sob}}(f)$, satisfies 
	\begin{equation}
		\var{I_N^{\mathrm{sob}}} \lesssim
		\begin{cases}
			N^{-1-\frac{1}{s} }  \left( \log N\right)^{s-1}, & \text{if } A^{*} < 0,\\	
			N^{-1 + \frac{(1-2A^{*})(-1+2A^{*}+\delta)}{s-2A^{*}-\delta}} \left( \log N \right)^{\frac{(s-1)(1-2A^{*} -\delta)}{s-2A^{*} -\delta }}, & \text{else},
		\end{cases}
	\end{equation}
	for $\delta > 0$. 
\end{Lemma}

\subsection{Partially axis-parallel sets}
This section considers the case when $\Omega = \Omega_{\mathfrak{u}} \times [\bm{a}_{-\mathfrak{u}}, \bm{b}_{-\mathfrak{u}})$, with $\mathfrak{u} \subsetneq 1:s$, is a partially axis-parallel set, where $\Omega_{\mathfrak{u}}$ admits a $\abs{\mathfrak{u}} - 1 $ dimensional Minkowski content. We summarize the results in the following lemma. 
\begin{Lemma}
	\label{lemma:sobol_partial_axis_parallel_variance_decay}
	For $f = g \cdot \mathbbm{1}_{\Omega}$, where $g$ satisfies the boundary growth condition~\eqref{eq:boundary_growth_condition},  $\Omega = \Omega_{\mathfrak{u}} \times [\bm{a}_{-\mathfrak{u}}, \bm{b}_{-\mathfrak{u}}) \subset [0, 1]^s$, and $\Omega_{\mathfrak{u}}$ admits an $\abs{\mathfrak{u}}-1$ dimensional Minkowski content, the variance of the RQMC estimator with scrambled Sobol' sequence $I_N^{\mathrm{sob}} \coloneqq I_N^{\mathrm{sob}}(f)$ satisfies 
	\begin{equation}
		\var{I_N^{\mathrm{sob}}} \lesssim
		\begin{cases} N^{-1 - \frac{1}{\abs{\mathfrak{u}} }} (\log N)^{\frac{s - 1}{\abs{\mathfrak{u}}}}, &\ \mathrm{if} \ A^{*} < 0,\\
			N^{-1 + \frac{(1-2A^{*})(-1+2A^{*}+\delta)}{\abs{\mathfrak{u}} - 2A^{*} -\delta}} (\log N)^{\frac{(s - 1)(1-2A^{*} -\delta)}{\abs{\mathfrak{u}}-2A^{*} - \delta}}, &\ \mathrm{otherwise},
		\end{cases}
	\end{equation}
	for $\delta > 0$. 
\end{Lemma}

The proof largely follows from the derivations in the previous section, and we defer it to Appendix~\ref{section:proof_sobol_partial_axis_parallel}. 
Lemma~\ref{lemma:sobol_partial_axis_parallel_variance_decay} improves the convergence rates compared to the results in~\cite{he2015convergence,he2018quasi_discontinuous_boundary} for both bounded and unbounded discontinuous integrands. 

\section{Numerical Examples}
\label{sec:numex}
Following Section~\ref{sec:discontinuous_integrand}, we consider the integrand of the form $f(\bm{t}) = \mathbbm{1}_{\bm{t} \in \Omega} g(\bm{t})$, where $g$ is derived from the numerical example in~\cite{ouyang2023quasi}:
\begin{equation}
	\label{eq:integrand}
	g(\bm{t}) = C \exp \left( M \norm{\Phi^{-1}(\bm{t})}_2^2 \right),
\end{equation}
where $M < 1/4$, and the normalizing constant $C = (1 - 2M)^{s/2}$. We consider various types of $\Omega$ in the following four examples and compare the observed convergence rates with theoretical predictions. In this section, we employ the lattice rule from~\cite{cools2006constructing} and Sobol' sequence from~\cite{joe2008constructing}. 

\begin{figure}[htbp]
	\centering
	\includegraphics[width=0.32\textwidth]{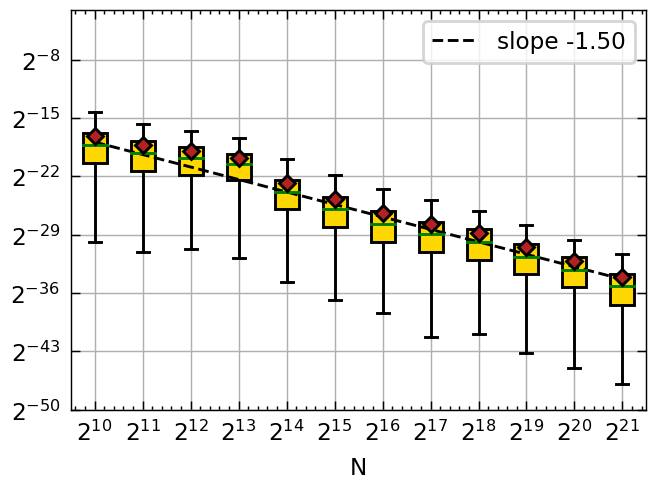}
	\includegraphics[width=0.32\textwidth]{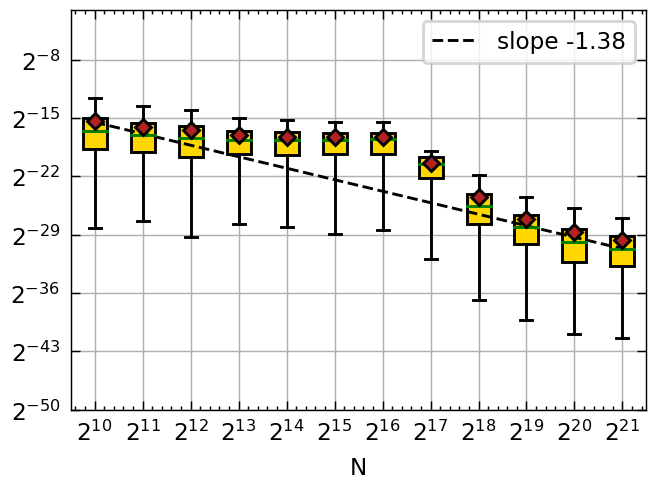}
	\includegraphics[width=0.32\textwidth]{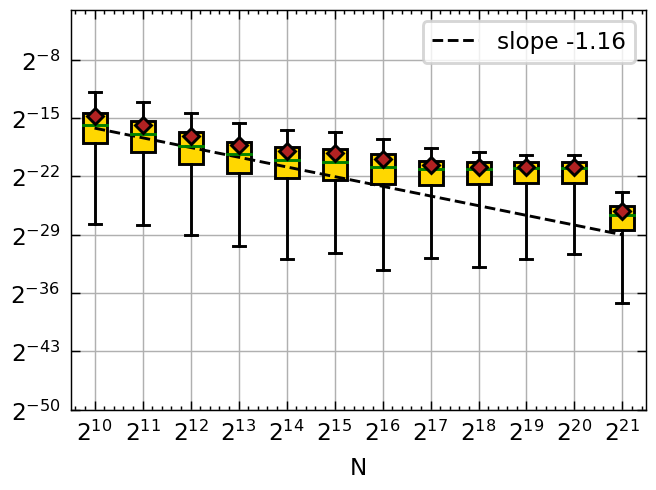}
	\hfill\\
	\includegraphics[width=0.32\textwidth]{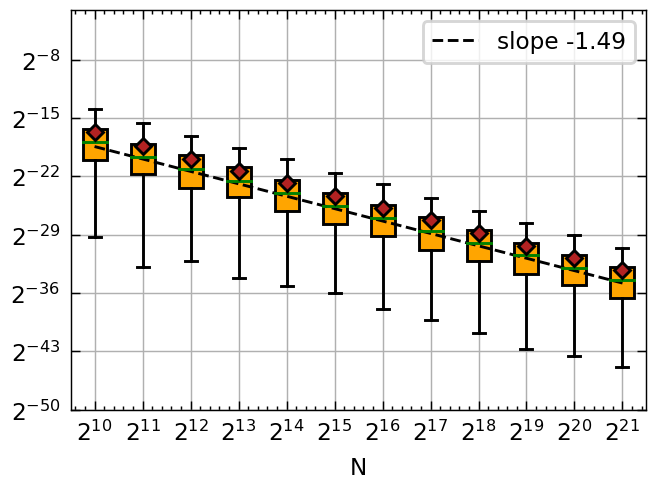}
	\includegraphics[width=0.32\textwidth]{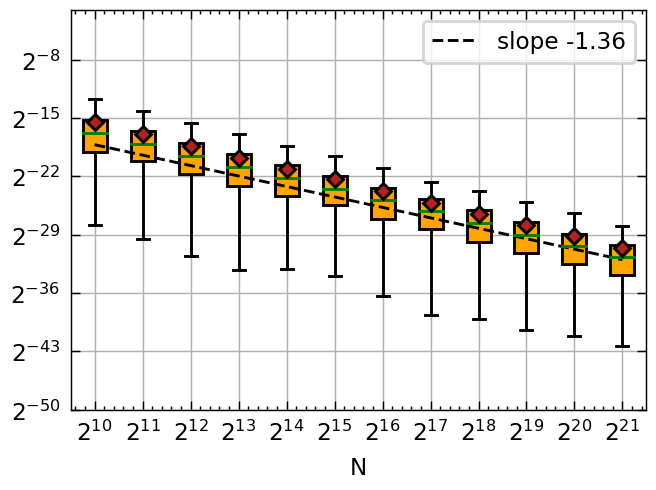}
	\includegraphics[width=0.32\textwidth]{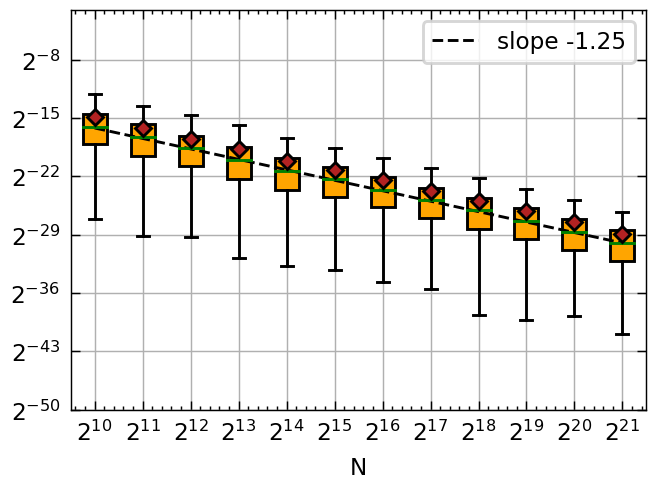}
	\hfill\\
	\caption{Ex1: Convergence of RQMC estimator variances. Top: RSLR; bottom: Sobol' sequence. The columns from left to right correspond to $s = 2, 3, 5$ respectively. {$M = 0.1$}. Each boxplot displays the 1st to 99th percentile of 2,048 samples. }
	\label{fig:ex2_s235_M1e-1}
\end{figure}

\begin{figure}[htbp]
	\centering
	\includegraphics[width=0.32\textwidth]{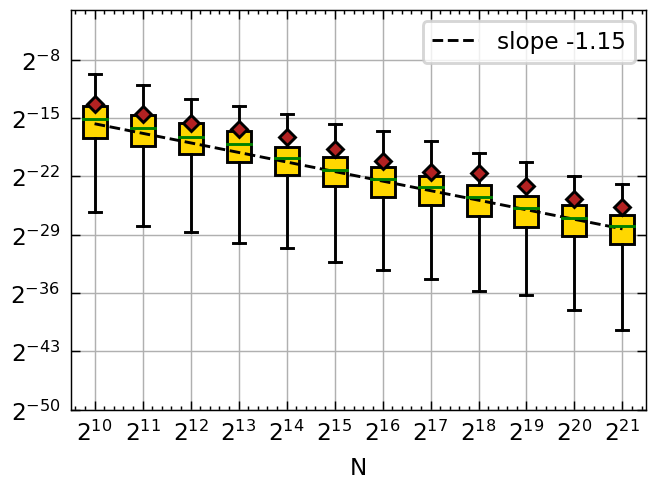}
	\includegraphics[width=0.32\textwidth]{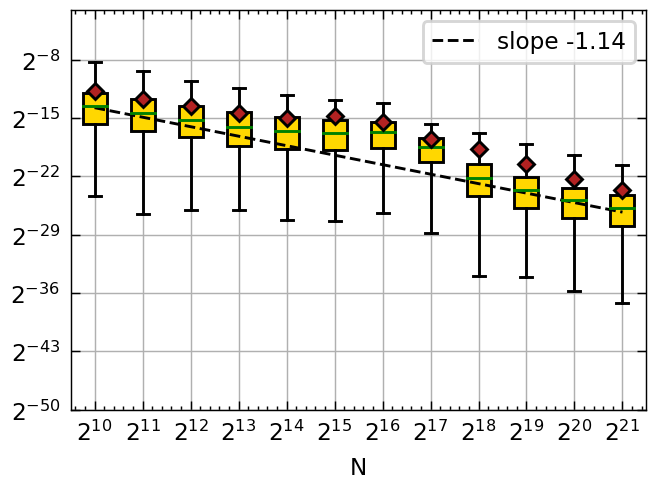}
	\includegraphics[width=0.32\textwidth]{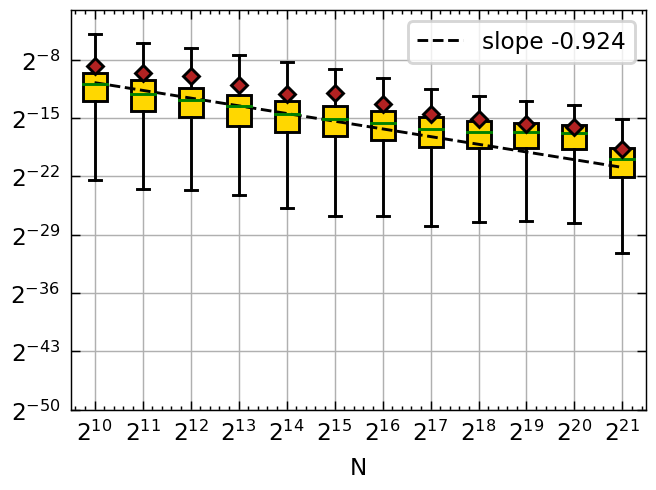}
	\hfill\\
	\includegraphics[width=0.32\textwidth]{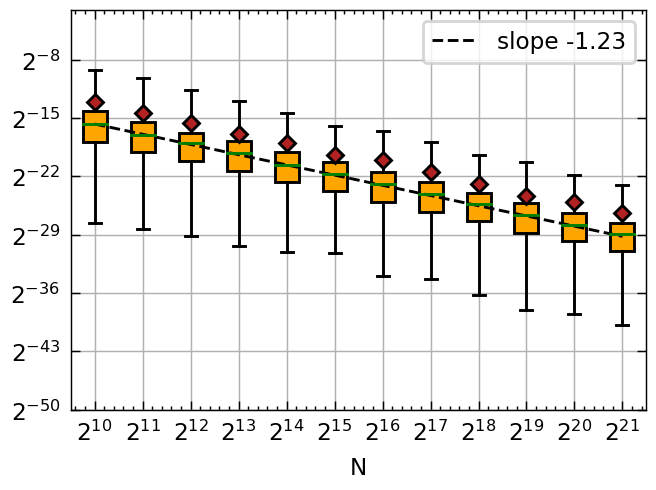}
	\includegraphics[width=0.32\textwidth]{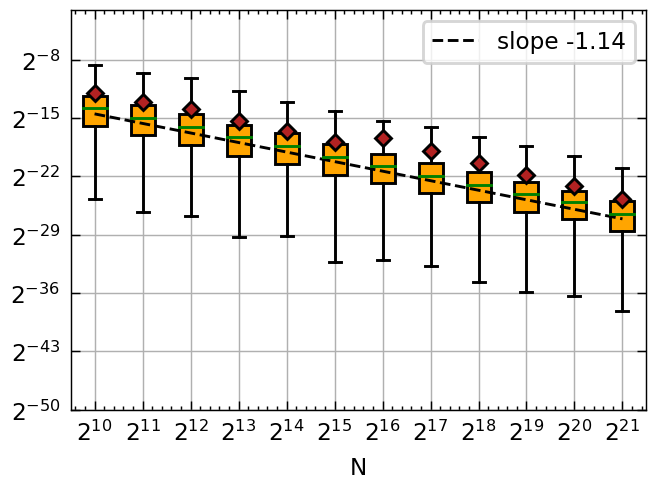}
	\includegraphics[width=0.32\textwidth]{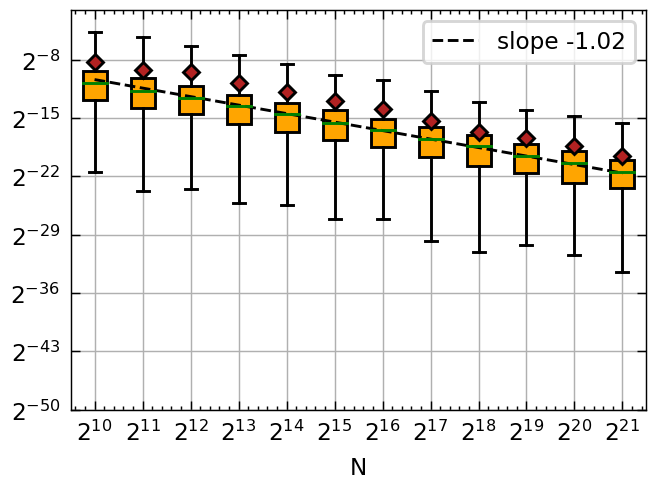}
	\hfill\\
	\caption{Ex1: Convergence of RQMC estimator variances. Top: RSLR; bottom: Sobol' sequence. The columns from left to right correspond to $s = 2, 3, 5$, with {$M = 0.2$} when $s = 2, 3$, and {$M = 0.225$} when $s = 5$. Each boxplot displays the 1st to 99th percentile of 2,048 samples. }
	\label{fig:ex2_s235_M2e-1}
\end{figure}

In Example 1, we analyze a hypersphere-type discontinuity with $\Omega = \{\bm{t} \in [0, 1]^s: \sum_{j=1}^s {t}_j^2 \leq 1 \}$. Figure~\ref{fig:ex2_s235_M1e-1} and~\ref{fig:ex2_s235_M2e-1} display boxplots for 2,048 samples of squared errors, $(I_N^{\mathrm{lat}} - I)^2$ and $(I_N^{\mathrm{sob}} - I)^2$, across various settings of $M$. 

In Figure~\ref{fig:ex2_s235_M1e-1}, when $M = 0.1$, both RSLR and the randomly scrambled Sobol' sequence are affected by the hypersphere-type discontinuity rather than the boundary growth conditions. For both sequences, the fitted convergence rates are close to $\mathcal{O}(N^{-1-\frac{1}{s}+\epsilon})$, faster than our theoretical bounds. When the value of $M$ is increased to $0.2$ for $s = 2, 3$ and $0.225$ for $s = 5$ in Figure~\ref{fig:ex2_s235_M2e-1}, the convergence rates are close to $\mathcal{O}(N^{-2+4M+\epsilon})$, suggesting the effects from the boundary growth condition, again exceeding our theoretical error bounds. The suboptimal observed rates for RSLR are due to the nonasymptotic regimes. Additionally, the RSLR presents large errors for certain $N$, suggesting that the generating vector $\bm{z}$ may be suboptimal for these cases. For instance, the reader may refer to~\cite{l2020tool} for optimizing the generating vector for RSLR.

\begin{figure}[htbp]
	\centering
	\includegraphics[width=0.32\textwidth]{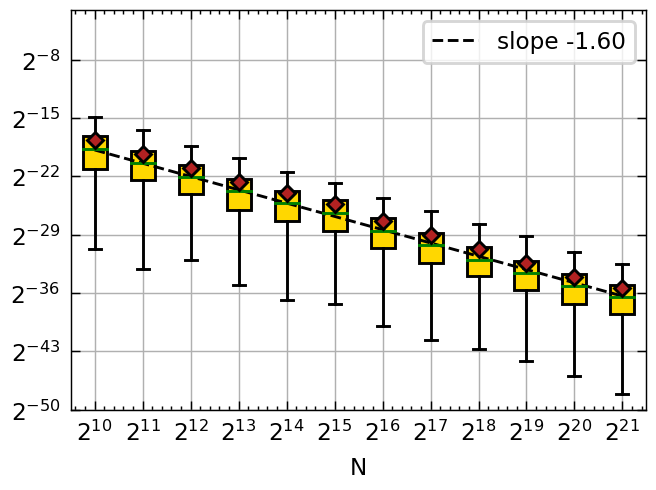}
	\includegraphics[width=0.32\textwidth]{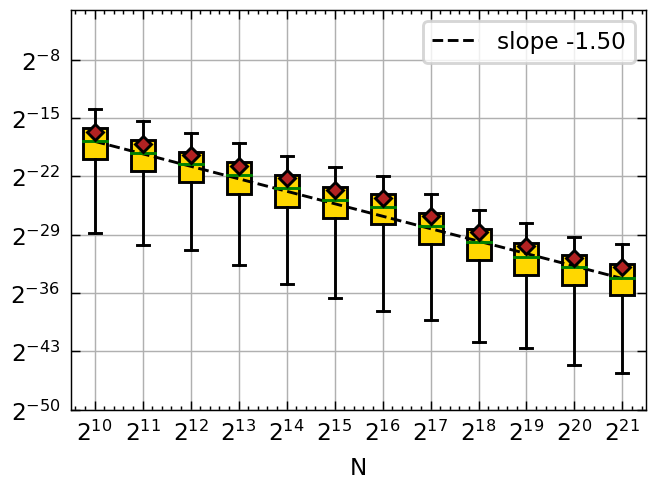}
	\includegraphics[width=0.32\textwidth]{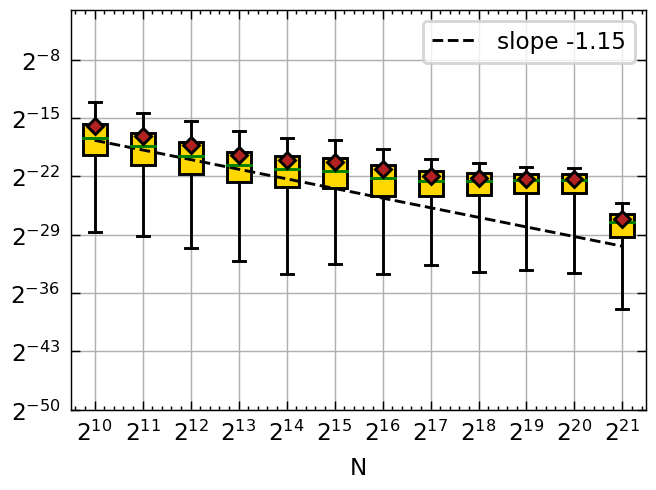}
	\hfill\\
	\includegraphics[width=0.32\textwidth]{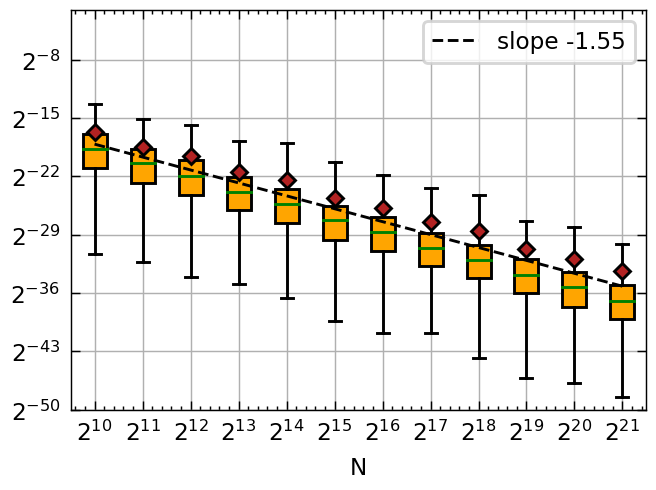}
	\includegraphics[width=0.32\textwidth]{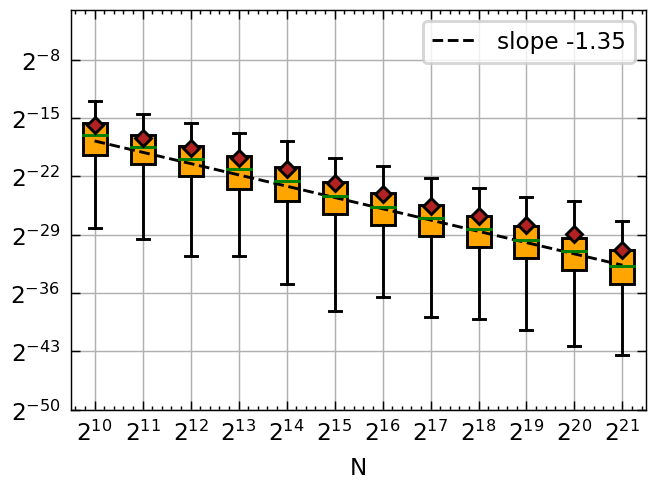}
	\includegraphics[width=0.32\textwidth]{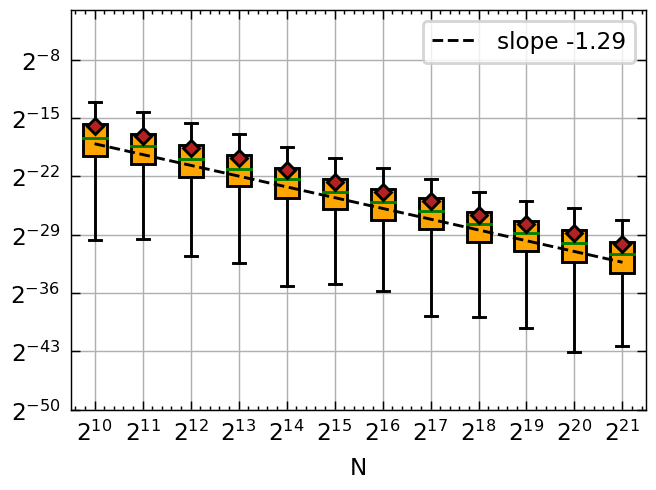}
	\hfill\\
	\caption{Ex2: Convergence of RQMC estimator variances. Top: RSLR; bottom: Sobol' sequence. The columns from left to right correspond to $s = 2, 3, 5$, respectively. {$M = 0.1$}. Each boxplot displays the 1st to 99th percentile of 2,048 samples.}
	\label{fig:ex1_s235_M1e-1}
\end{figure}

\begin{figure}[htbp]
	\centering
	\includegraphics[width=0.32\textwidth]{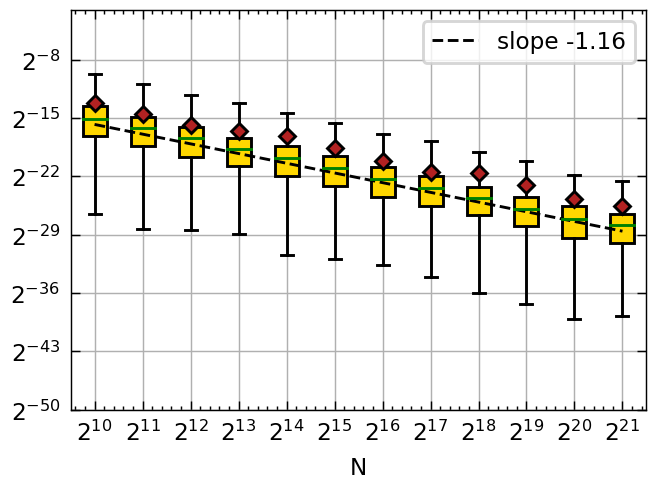}
	\includegraphics[width=0.32\textwidth]{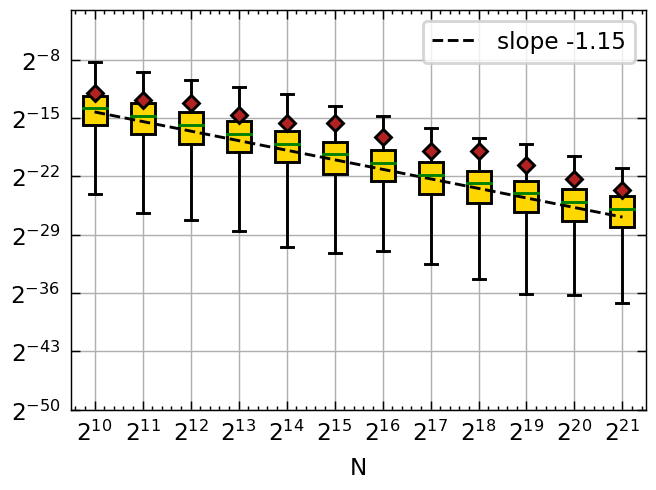}
	\includegraphics[width=0.32\textwidth]{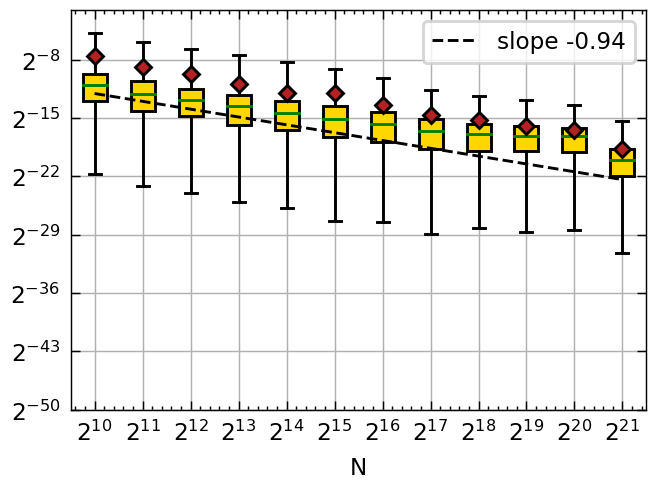}
	\hfill\\
	\includegraphics[width=0.32\textwidth]{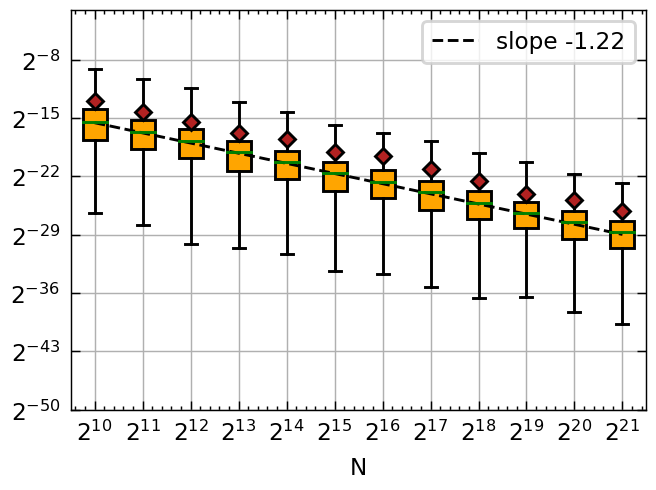}
	\includegraphics[width=0.32\textwidth]{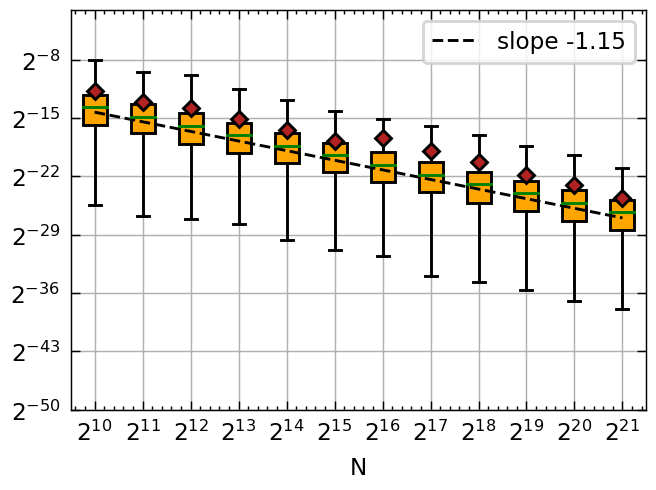}
	\includegraphics[width=0.32\textwidth]{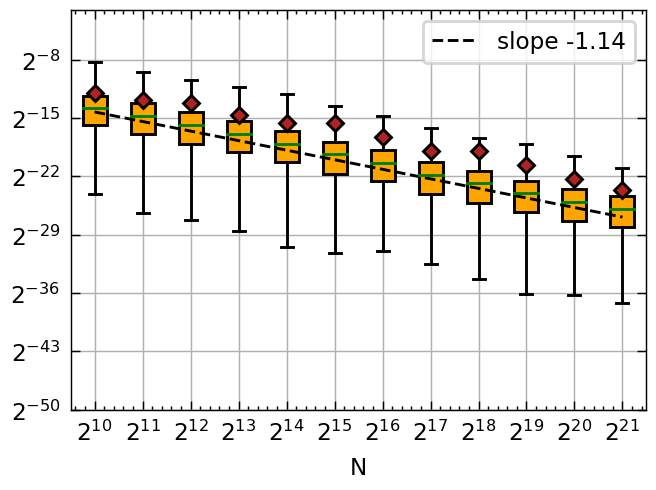}
	\hfill\\
	\caption{Ex2: Convergence of RQMC estimator variances. Top: RSLR; bottom: Sobol' sequence. The columns from left to right correspond to $s = 2, 3, 5$, with {$M = 0.2$} when $s = 2, 3$, and {$M = 0.225$} when $s = 5$. Each boxplot displays the 1st to 99th percentile of 2,048 samples.} 
	\label{fig:ex1_s235_M2e-1}
\end{figure}

In Example 2, we consider the hyperplane-type discontinuity with $\Omega = \{\bm{t} \in [0, 1]^s: \sum_{j=1}^s {t}_j \leq 1 \}$. All other settings remain the same as in Example 1. Boxplots of squared errors are shown in Figure~\ref{fig:ex1_s235_M1e-1} and~\ref{fig:ex1_s235_M2e-1}. 
The Sobol' sequences behave similarly to those in Example 1 with the hypersphere-type discontinuity. However, the behavior of RSLR is worth discussing. When $M = 0.1$ and $s=2$, the fitted convergence $-1.60$ closely aligns with the convergence rate $\mathcal{O}(N^{-2+4M+\epsilon})$, suggesting that the boundary growth condition, rather than the discontinuity, dominates the convergence rate --- a contrast to Example 1. Moreover, when $M = 0.1$ and $s=3, 5$, the fitted convergence rates are close to $\mathcal{O}(N^{-1-\frac{1}{s - 1}+\epsilon})$, as if the dimension is reduced by 1 in the discontinuity effect. This difference from Example 1 arises due to the distinct nature of the discontinuity boundaries. When $M$ is increased as shown in Figure~\ref{fig:ex1_s235_M2e-1}, the boundary growth conditions dominate the convergence rates, resembling the behavior in Example 1. 

\begin{figure}[htbp]
	\centering
	\includegraphics[width=0.32\textwidth]{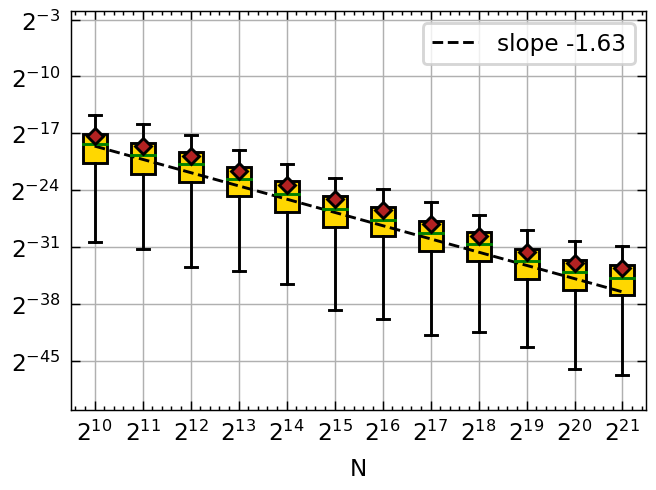}
	\includegraphics[width=0.32\textwidth]{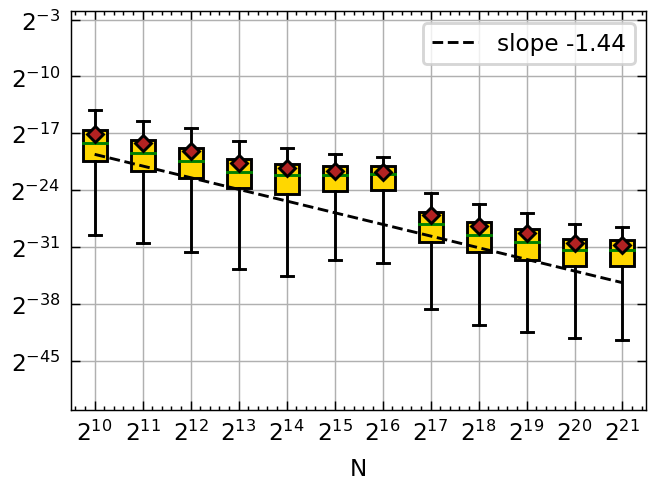}
	\includegraphics[width=0.32\textwidth]{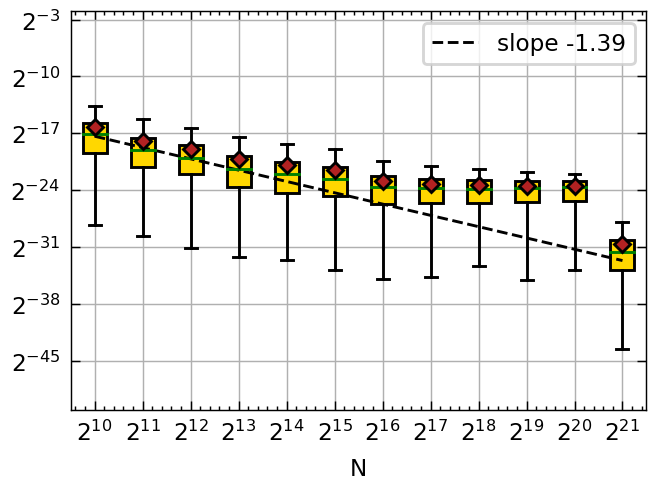}
	\hfill\\
	\includegraphics[width=0.32\textwidth]{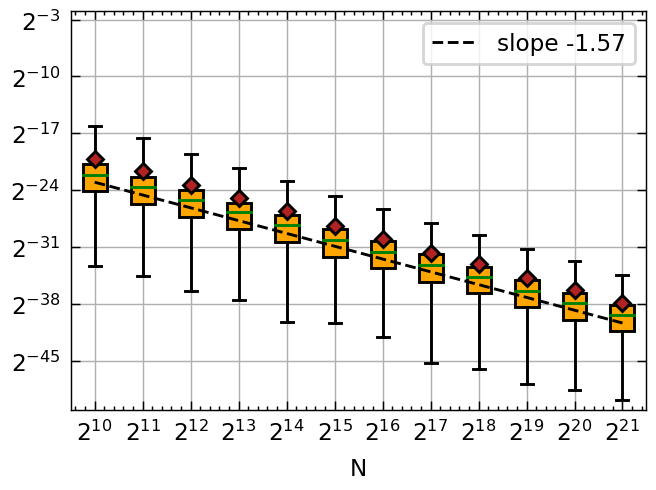}
	\includegraphics[width=0.32\textwidth]{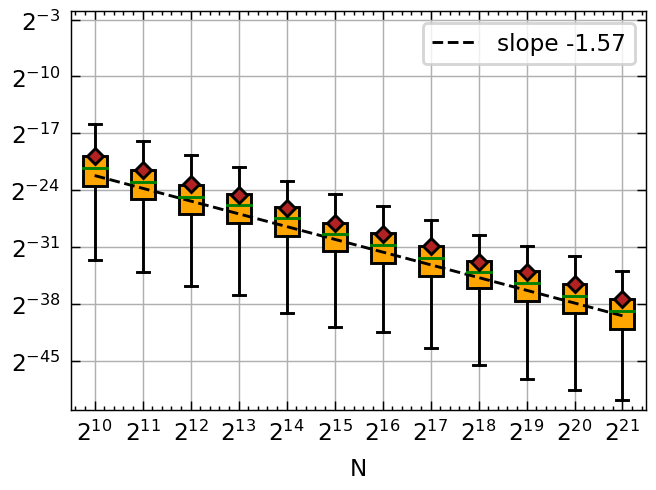}
	\includegraphics[width=0.32\textwidth]{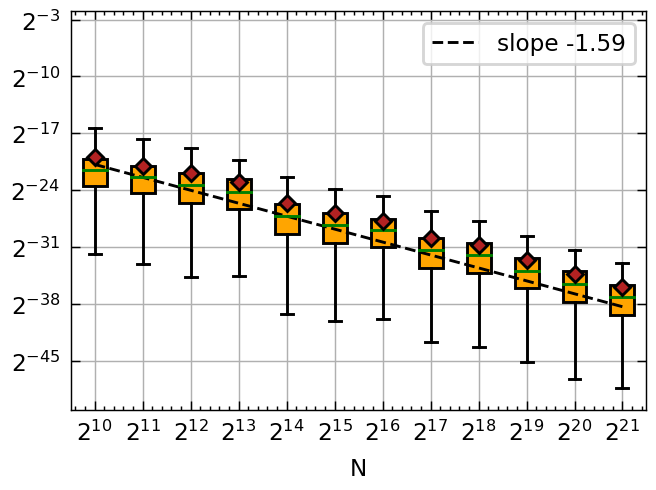}
	\hfill\\
	\caption{Ex3: Convergence of RQMC estimator variances. Top: RSLR, bottom: Sobol' sequence. The columns from left to right correspond to $s = 2, 3, 5$, respectively. {$M = 0.1$}. Each boxplot displays the 1st to 99th percentile of 2,048 samples.}
	\label{fig:ex3_s235_M1e-1}
\end{figure}

In Example 3, we consider an axis-parallel discontinuity. Specicially, we consider $\Omega = \{ \bm{t} \in [0, 1]^s: \min_{j=1}^s t_j \geq 1/2 \}$ and let $M = 0.1$. Boxplots of squared errors are shown in Figure~\ref{fig:ex3_s235_M1e-1}. 

From Figure~\ref{fig:ex3_s235_M1e-1}, the convergence rates of Sobol' sequences for all three dimension settings closely align with our theoretical rates $\mathcal{O}(N^{-2+4M+\epsilon})$, indicating the boundary growth condition dominates the convergence rate. More nonasymptotic effects appear in the RSLR when dimension $s$ increases. 

\begin{figure}[htbp]
	\centering
	\includegraphics[width=0.32\textwidth]{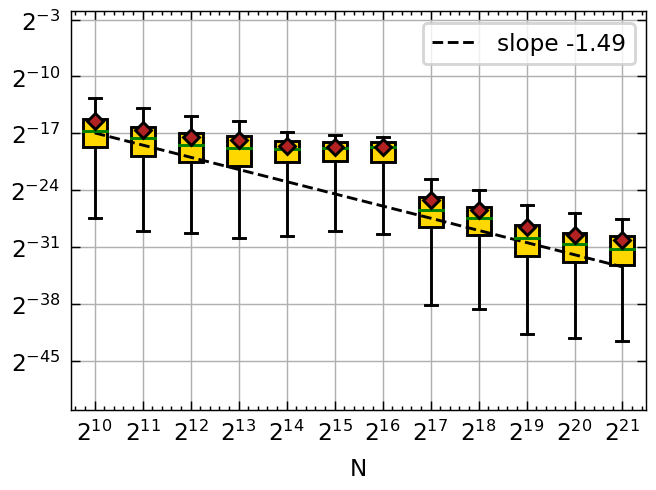}
	\includegraphics[width=0.32\textwidth]{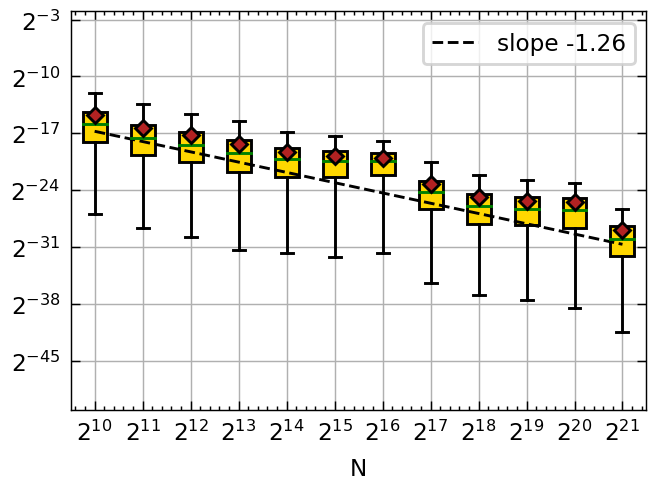}
	\includegraphics[width=0.32\textwidth]{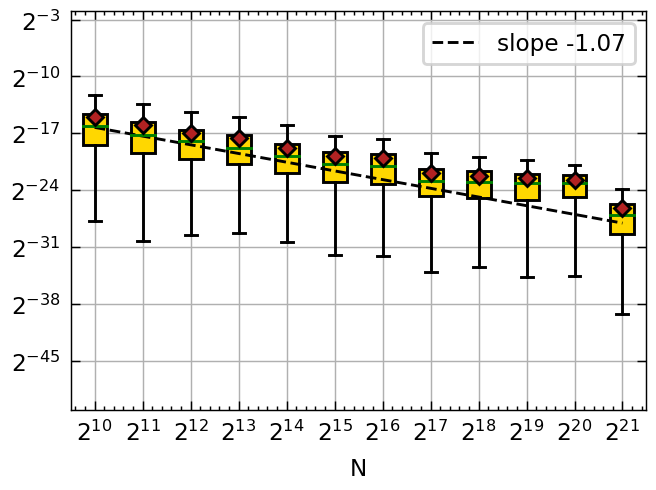}
	\hfill\\
	\includegraphics[width=0.32\textwidth]{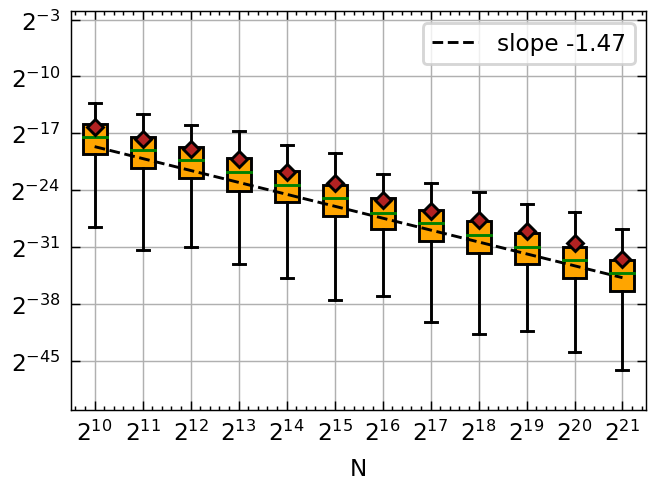}
	\includegraphics[width=0.32\textwidth]{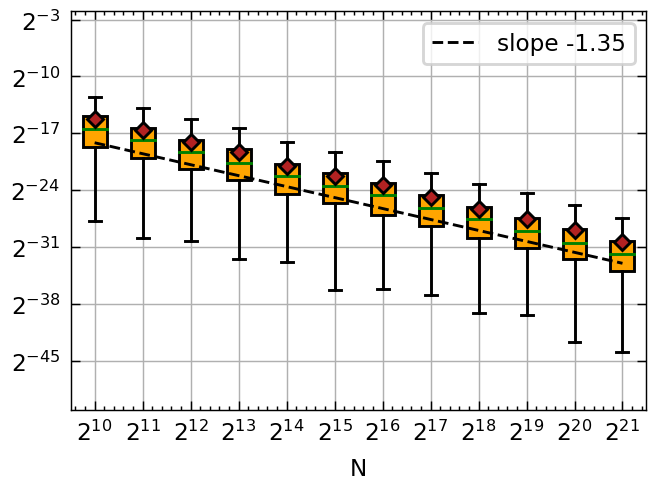}
	\includegraphics[width=0.32\textwidth]{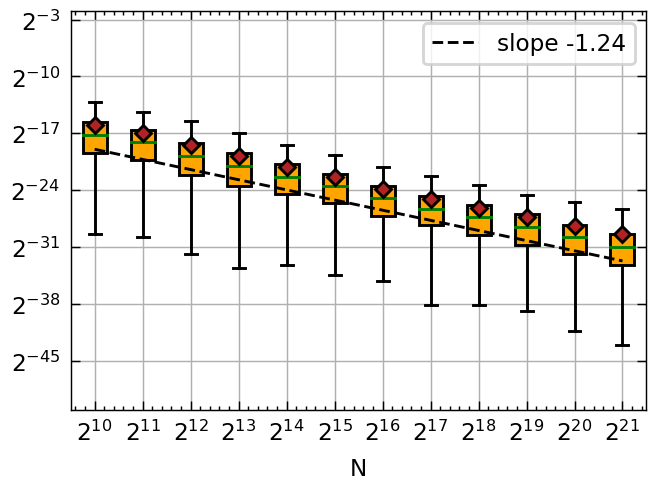}
	\hfill\\
	\caption{Ex4: Convergence of RQMC estimator variances. Top: RSLR, bottom: Sobol' sequence. The columns from left to right corresponds to {$s = 4, 5, 7$}, respectively. {$M = 0.1$}. Each boxplot displays the 1st to 99th percentile of 2,048 samples.} 
	\label{fig:ex4_s235_M1e-1}
\end{figure}

In Example 4, we consider a partial axis-parallel discontinuity. We consider $\Omega = \{ \bm{t} \in [0, 1]^s: \min_{j=1}^2 {t}_j \geq 1/2, \sum_{i=3}^s {t}_i^2 \geq 1 \}$ with $M = 0.1$ and dimensions $s = 4, 5, 7$. The boxplots of squared errors are shown in Figure~\ref{fig:ex4_s235_M1e-1}. 
The convergence rates from Figure~\ref{fig:ex4_s235_M1e-1} are similar to those in Figure~\ref{fig:ex1_s235_M1e-1} from Example 1. This indicates that in the case of partial axis-parallel discontinuity, the effects of the axis-parallel dimensions are marginal, consistent with our theoretical results.

This section presents numerical examples for unbounded and discontinuous integrands. The numerical results for either unbounded or discontinuous integrands can be referred to~\cite{he2015convergence,ouyang2023quasi}, which align with the theoretical results in this work. 

\section{Conclusions}
\label{sec:conclusion}
This work presents a comparative study of the RSLR and randomly scrambled Sobol' sequences through spectral analysis. We analyze integrands with infinite Hardy--Krause variation due to the boundary unboundedness or interior discontinuities. We establish sufficient conditions under which these integrands achieve specific convergence rates and compare our findings to classical results. Using Owen's boundary growth conditions, we provide fractional-order convergence rates that bridge the gaps between integer differentiability orders and their associated convergence rates. 

Both low-discrepancy sequences under consideration are insensitive to {axis-parallel} discontinuities. For general discontinuous integrands with Sobol' sequences, the convergence rates observed in the numerical results are faster than our theoretical bounds. One potential reason could be that the H\"older inequality applied in~\eqref{eq:holder_inequality} does not provide optimal upper bounds. Another reason can be the loose upper bound in Lemma~\ref{lemma:sobol_axis_parallel_multiple}. Future work could aim to refine these upper bounds. 

The effect of the discontinuity interface on the RSLR remains unclear. \cite{greenblatt2021fourier} studies the Fourier transform of indicator functions on $\mathbb{R}^s$ and reveals its dependency on the boundary curvature. Future work could extend the analysis in~\cite{greenblatt2021fourier} to study the decay of Fourier coefficients, for instance, for the unbounded and discontinuous integrands examined in this work. 

Furthermore, this study could facilitate the variance reduction techniques in the RQMC settings. For instance, in the context of importance sampling for QMC---as discussed in~\cite{he2022error, liu2023nonasymptotic}---future efforts could consider leveraging spectral properties in additional to minimizing the second-order moment when searching for the optimal proposal distribution. 
\section*{Acknowledgments}
This work is funded by the Alexander von Humboldt Foundation and King Abdullah University of Science and Technology (KAUST) Office of Sponsored Research (OSR) under Award No. OSR-2019-CRG8-4033. 
This work utilized the resources of the Supercomputing Laboratory at King Abdullah University of Science and Technology (KAUST) in Thuwal, Saudi Arabia. The author extends gratitude to Prof. Raúl Tempone for his helpful discussions and to Prof. Zhijian He for identifying errors in Lemma 4.5 of the original draft. The author also thanks the two anonymous reviewers for their constructive comments, which greatly contributed to the improvement of this manuscript.

\bibliographystyle{siam}
\bibliography{sampling, bibliography_QMC_theory, bibliography_QMC_finance}
\include{appendix.tex}

\end{document}

%% file: commands_blue_noise.tex
\usepackage{mathtools}

\providecommand{\E}[1]{{\ensuremath{\mathbb{E}}\mspace{-2mu}\left[#1\right]}}    
\providecommand{\var}[1]{{\ensuremath{\mathrm{Var}}\mspace{-2mu}\left[#1\right]}}

\providecommand{\wal}{\ensuremath{{\mathrm{wal}}}} 

\DeclarePairedDelimiter\abs{\lvert}{\rvert}%
\DeclarePairedDelimiter\norm{\lVert}{\rVert}%

\newtheorem{Remark}{Remark}
\newtheorem{Definition}{Definition}

\newtheorem{Assumption}{Assumption}
\newtheorem{Theorem}{Theorem}

\newtheorem{Lemma}{Lemma}

\newcommand{\RN}[1]{%
	\textup{\uppercase\expandafter{\romannumeral#1}}%
}

%% file: appendix.tex
\appendix
\numberwithin{equation}{section}

\section{Proof to Lemma~\ref{lemma:decay_cn_1d}}
\label{section:proof_decay_cn_1d}
\begin{proof}
	In 1-D case, the Fourier coefficients are given by
   \begin{equation}
	   \begin{split}
		   c_n &= \int_{\mathbb{T}} f(t) e^{-i2\pi n t} dt = \int_{\mathbb{T}} f(t + \frac{1}{2n}) e^{-i2\pi n (t + \frac{1}{2n})} dt  = \int_{\mathbb{T}} -f(t + \frac{1}{2n}) e^{-i2\pi n t } dt.
	   \end{split}
   \end{equation}
   We can derive the following bound for $\abs{c_n}$:
   \begin{equation}
	\label{eq:bound_cn_1d}
	   \begin{split}
		   \abs{c_n} &= \frac{1}{2} \abs*{\int_{\mathbb{T}} f(t) e^{-i2\pi n t} dt + \int_{\mathbb{T}} -f(t + \frac{1}{2n}) e^{-i2\pi n t } dt} \leq \frac{1}{2} \int_{\mathbb{T}} \abs*{ f(t) - f(t + \frac{1}{2n}) }dt\\
		   &\leq \frac{1}{2} \int_{\mathbb{T}} \left( \mathbbm{1}_{t \in [0, 1 - \frac{1}{2n})} \int_{t}^{t+\frac{1}{2n}} \abs*{\frac{\partial f}{\partial t_0} }dt_0 + \mathbbm{1}_{t \in [1 - \frac{1}{2n}, 1)} \int_{t+\frac{1}{2n} - 1}^{t} \abs*{\frac{\partial f}{\partial t_0} }dt_0 \right)dt .
	   \end{split}
   \end{equation}
   Next, we apply the upper bound $\abs*{\frac{\partial f}{\partial t_0} } \lesssim \phi(t_0)^{-A-1}$, and consider the following four situations for $A \neq 0$ :(i) $t \in [0, \frac{1}{2} - \frac{1}{2n})$
   \begin{equation}
	   \begin{split}
		 \int_{0}^{\frac{1}{2} - \frac{1}{2n}}  \int_{t}^{t + \frac{1}{2n}} \abs*{\frac{\partial f}{\partial t_0}} dt_0 dt &\lesssim \int_{0}^{\frac{1}{2} - \frac{1}{2n}}  \int_{t}^{t + \frac{1}{2n}} t_0^{-A-1} dt_0 dt
		   = \frac{1}{A(1-A)} \left[ \left( \frac{1}{2} - \frac{1}{2n} \right)^{1-A} - 2^{A-1} + \left(\frac{1}{2n}\right)^{1-A} \right]. 
	   \end{split}
	   \label{eq:bound_cn_1d_1}
   \end{equation}
   When (ii) $t \in [\frac{1}{2} - \frac{1}{2n}, \frac{1}{2})$, we have
   \begin{equation}
	   \begin{split}
		\int_{\frac{1}{2} - \frac{1}{2n}}^{\frac{1}{2}} \int_{t}^{t + \frac{1}{2n}} \abs*{\frac{\partial f}{\partial t_0}} dt_0 dt &= \int_{\frac{1}{2} - \frac{1}{2n}}^{\frac{1}{2}} \left(  \int_{t}^{\frac{1}{2}} \abs*{\frac{\partial f}{\partial t_0}} dt_0 + \int_{\frac{1}{2}}^{t + \frac{1}{2n}}  \abs*{\frac{\partial f}{\partial t_0}} dt_0 \right) dt \\
		   &\lesssim  -\frac{1}{A}\cdot \frac{2^A}{n} +  \frac{2}{A(1-A)} \left[ 2^{A-1} - 2\cdot\left( \frac{1}{2} - \frac{1}{2n} \right)^{1-A}\right].
	   \end{split}
   \end{equation}
   Similar to case (i), when (iii) $t \in [\frac{1}{2}, 1 - \frac{1}{2n})$, we have
   \begin{equation}
	   \int_{\frac{1}{2}}^{1 - \frac{1}{2n}} \int_{t}^{t + \frac{1}{2n}} \abs*{\frac{\partial f}{\partial t_0}} dt_0 dt \lesssim \frac{1}{A(1-A)} \left[ \left( \frac{1}{2} - \frac{1}{2n} \right)^{1-A} - 2^{A-1} + \left(\frac{1}{2n}\right)^{1-A} \right]. 
   \end{equation}
   In case (iv) $t \in [1 - \frac{1}{2n}, 1)$, we have
   \begin{equation}
	   \begin{split}
		\int_{1 - \frac{1}{2n}}^{1} \int_{t + \frac{1}{2n} - 1}^{\frac{1}{2}}  \abs*{\frac{\partial f}{\partial t_0}} dt_0 + \int_{\frac{1}{2}}^{t} \abs*{\frac{\partial f}{\partial t_0}} dt_0 dt
		&\lesssim -\frac{1}{A}\cdot \frac{2^A}{2n} +  \frac{2}{A(1-A)}  \left( \frac{1}{2n} \right)^{1-A}, 
	   \end{split}
	   \label{eq:bound_cn_1d_iv_2}
   \end{equation}
   where in Equation~\eqref{eq:bound_cn_1d_iv_2} we have used the fact that $f$ is $1$-periodic.
   Combining~\eqref{eq:bound_cn_1d_1}-\eqref{eq:bound_cn_1d_iv_2}, we have
   \begin{equation}
	   \abs{c_n} \lesssim \frac{4}{A(1-A)} \left( \frac{1}{2n} \right)^{1-A} - \frac{2^{A+1}}{An} \lesssim \begin{cases}
			n^{-1+A} & \text{if } A > 0,\\
		   n^{-1} & \text{if } A < 0.
	   \end{cases}  
   \end{equation}
   Similarly, in the case $A = 0$, we consider the following situations:(i) $t \in [0, \frac{1}{2} - \frac{1}{2n})$
   \begin{equation}
	   \begin{split}
		   \int_{0}^{\frac{1}{2} - \frac{1}{2n}}  \int_{t}^{t + \frac{1}{2n}} \abs*{\frac{\partial f}{\partial t_0}} dt_0 dt &\lesssim \int_{0}^{\frac{1}{2} - \frac{1}{2n}}  \int_{t}^{t + \frac{1}{2n}} t_0^{-1} dt_0 dt = \int_{0}^{\frac{1}{2} - \frac{1}{2n}} \log\left( t + \frac{1}{2n} \right) - \log{t} dt\\
		   &= \frac{1}{2}\log{\frac{1}{2}} - \frac{1}{2} - \left( \frac{1}{2n}\log{\frac{1}{2n}} - \frac{1}{2n} \right) - \left( \frac{1}{2} - \frac{1}{2n} \right) \log{\left(\frac{1}{2} - \frac{1}{2n}\right)}  + \frac{1}{2} - \frac{1}{2n}. 
	   \end{split}
	   \label{eq:bound_cn_1d_1_A0}
   \end{equation}
   In case (ii) when $t \in [\frac{1}{2} - \frac{1}{2n}, \frac{1}{2})$, we have
   \begin{equation}
	   \begin{split}
		   \int_{\frac{1}{2} - \frac{1}{2n}}^{\frac{1}{2}} \int_{t}^{t + \frac{1}{2n}} \abs*{\frac{\partial f}{\partial t_0}} dt_0 dt & \leq\int_{\frac{1}{2} - \frac{1}{2n}}^{\frac{1}{2}} \int_{t}^{\frac{1}{2}} \abs*{\frac{\partial f}{\partial t_0}} dt_0 + \int_{\frac{1}{2}}^{t + \frac{1}{2n}} \abs*{\frac{\partial f}{\partial t_0}} dt_0 dt \\
		   &\lesssim   \frac{1}{n} \log{\frac{1}{2}} - 2 \left[ \frac{1}{2} \log{\frac{1}{2}} - \frac{1}{2} - \left(\frac{1}{2} - \frac{1}{2n}\right) \log{\left( \frac{1}{2} - \frac{1}{2n}\right)} + \frac{1}{2} - \frac{1}{2n} \right].
	   \end{split}
   \end{equation}
   The case (iii) when $t \in [\frac{1}{2}, 1 - \frac{1}{2n})$ is similar to case (i), and we have
   \begin{equation}
	\begin{split}
		\int_{0}^{\frac{1}{2} - \frac{1}{2n}}  \int_{t}^{t + \frac{1}{2n}} \abs*{\frac{\partial f}{\partial t_0}} dt_0 dt &\lesssim \frac{1}{2}\log{\frac{1}{2}} - \frac{1}{2} - \left( \frac{1}{2n}\log{\frac{1}{2n}} - \frac{1}{2n} \right) - \left( \frac{1}{2} - \frac{1}{2n} \right) \log{\left(\frac{1}{2} - \frac{1}{2n}\right)}  + \frac{1}{2} - \frac{1}{2n}. 
	\end{split}
   \end{equation}
   In case (iv) when $t \in [1 - \frac{1}{2n}, 1)$, we have
   \begin{equation}
	   \begin{split}
		   \int_{1-\frac{1}{2n}}^1 \int_{t + \frac{1}{2n} - 1}^{t} \abs*{\frac{\partial f}{\partial t_0}} dt_0 dt &\leq \int_{1 - \frac{1}{2n}}^{1} \int_{t + \frac{1}{2n} - 1}^{\frac{1}{2}}  \abs*{\frac{\partial f}{\partial t_0}} dt_0 + \int_{\frac{1}{2}}^{t} \abs*{\frac{\partial f}{\partial t_0}} dt_0 dt \lesssim  \frac{1}{n} \log{\frac{1}{2}}.
	   \end{split}
	   \label{eq:bound_cn_1d_iv_A0}
   \end{equation}
   Combining~\eqref{eq:bound_cn_1d_1_A0}-\eqref{eq:bound_cn_1d_iv_A0}, we have in case $A = 0$,
   \begin{equation}
	   \abs{c_n} \lesssim -\left(\frac{1}{n}\log{\frac{1}{2n}} - \frac{1}{n} \right) + \frac{2}{n} \log{\frac{1}{2}} \lesssim n^{-1}\log n.
   \end{equation}
\end{proof}

\section{Proof to Lemma~\ref{lemma:decay_cn_multidimension_axis_parallel}}
\label{section:proof_decay_cn_multidimension_axis_parallel}
\begin{proof}

	For $\bm{t} \in [0, 1)^s$ and $\bm{n} \in \mathbb{N}^s$, we consider a nonempty set $\mathfrak{u} \subseteq 1:s$ such that $\bm{n}_{\mathfrak{u}} = \bm{n}^{\mathfrak{u}}:0^{-\mathfrak{u}}$. Without loss of generality, we assume that ${1} \in \mathfrak{u}$. 
	
	We proceed with a proof by induction. For a function $f = \mathbbm{1}_{\Omega} \cdot g$, where $g$ satisfies the boundary growth condition~\eqref{assumption:owen_boundary_growth_condition} with $0 < A_j < 1/2$ for $j = 1, \dotsc, s$ and the set $\Omega = [\bm{a}, \bm{b}] \subseteq [0, 1]^s$ is axis-parallel, we assume
	\begin{equation}
		\label{eq:bound_cn_multidimension_axis_parallel_induction}
		\int_{[0, 1]^s}  \abs*{\Delta(f; \bm{t}, \bm{t} + \frac{1}{2\bm{n}_{\mathfrak{u}}})} d\bm{t} \lesssim  \prod_{j=1}^s {n}_j^{-1+A_j}. 
	\end{equation}
	The base case for $s = 1$ is already established in Lemma~\ref{lemma:decay_cn_1d}. 
	
	Now, we assume the statement~\eqref{eq:bound_cn_multidimension_axis_parallel_induction} holds for $s-1$ dimensions. In the case of $s$ dimensions, we consider to split the integration domain into two parts, depending on whether the interval $[{t}_1, {t}_1 + \frac{1}{2{n}_{1}}]$ intersects with the boundary of $\Omega_1$ or the boundary of $[0, 1]$, where the function $f$ becomes discontinuous. We have the following decomposition:
	\begin{equation}
		\label{eq:bound_cn_multidimension_axis_parallel}
		\begin{split}
			\int_{[0, 1]^s}  \abs*{\Delta(f; \bm{t}, \bm{t} + \frac{1}{2\bm{n}_{\mathfrak{u}}})} d\bm{t} &=  \int_{[0, 1]^{s-1}} \int_0^1 \mathbbm{1}_{[{t}_1, {t}_1 + \frac{1}{2 {n}_{1}}] \cap \partial (\bar{\Omega}_{1}) \neq \varnothing} \abs*{\Delta(f; \bm{t}, \bm{t} + \frac{1}{2\bm{n}_{\mathfrak{u}}}) } d {t}_1 d\bm{t}_{-1}\\
			&+ \int_{[0, 1]^{s-1}} \int_0^1 \mathbbm{1}_{[{t}_1, {t}_1 + \frac{1}{2 {n}_{1}}] \cap \partial (\bar{\Omega}_{1}) = \varnothing} \abs*{\Delta(f; \bm{t}, \bm{t} + \frac{1}{2\bm{n}_{\mathfrak{u}}}) } d {t}_1 d\bm{t}_{-1}.
		\end{split}
	\end{equation}
	Notice that in~\eqref{eq:bound_cn_multidimension_axis_parallel}, we define $\partial (\bar{\Omega}_1) \coloneqq \partial (\Omega_1) \cup \partial [0, 1]$. We denote $\Omega_1$ the projection of $\Omega$ onto dimension 1 and $\partial (\Omega_1)$ the boundary of $\Omega_1$. Specifically, $\partial (\Omega_1) = \{a_1, b_1\}$ and $\partial (\bar{\Omega}_1) = \{0, a_1, b_1, 1\}$, for $0 < a_1 < b_1 < 1$. 
	
	In the following, we first consider the case when ${[{t}_1, {t}_1 + \frac{1}{2 {n}_{1}}] \cap \partial (\bar{\Omega}_{1}) \neq \varnothing}$. 
	We denote $f_{{t}_1} (\bm{t}_{-1}) \coloneqq f(\bm{t}_{-1}\mid {t}_1)$.
	Notice that we have the following decomposition in the direction of $\bm{t}_1$ of the alternating sum, 
We have,
\begin{equation}
	\begin{split}
		{\Delta(f; \bm{t}, \bm{t} + \frac{1}{2\bm{n}_{\mathfrak{u}}}) } 
		& = {\Delta(f_{t_1} - f_{t_1 + \frac{1}{2n_1}}; \bm{t}_{-1}, \bm{t}_{-1} + \frac{1}{2\bm{n}_{\mathfrak{u} -1}}) },
	\end{split}
\end{equation}
and we proceed the first part of Equation~\eqref{eq:bound_cn_multidimension_axis_parallel} as
\begin{equation}
	\label{eq:bound_cn_multidimension_axis_parallel_case_1}
	\begin{split}
		& \int_{[0, 1]^{s-1}} \int_0^1 \mathbbm{1}_{[{t}_1, {t}_1 + \frac{1}{2 {n}_{1}}] \cap \partial (\bar{\Omega}_{1}) \neq \varnothing} \abs*{\Delta(f; \bm{t}, \bm{t} + \frac{1}{2\bm{n}_{\mathfrak{u}}}) } d {t}_1 d\bm{t}_{-1}\\
		&\leq \int_{[0, 1]^{s-1}} \int_{0}^{1} \mathbbm{1}_{[{t}_1, {t}_1 + \frac{1}{2 {n}_{1}}] \cap \partial (\bar{\Omega}_{1}) \neq \varnothing} \abs*{\Delta(f_{t_1}; \bm{t}_{-1}, \bm{t}_{-1} + \frac{1}{2\bm{n}_{\mathfrak{u} -1}}) } + \abs*{\Delta(f_{t_1 + \frac{1}{2n_1}}; \bm{t}_{-1}, \bm{t}_{-1} + \frac{1}{2\bm{n}_{\mathfrak{u} -1}}) } d{t}_1 d\bm{t}_{-1}\\
		&= \int_{[0, 1]^{s-1}} \int_{0}^{1} \mathbbm{1}_{[{t}_1, {t}_1 + \frac{1}{2 {n}_{1}}] \cap \partial (\bar{\Omega}_{1}) \neq \varnothing} \phi(t_1)^{-A_1} \abs*{\Delta(\frac{f_{t_1}}{\phi(t_1)^{-A_1}}; \bm{t}_{-1}, \bm{t}_{-1} + \frac{1}{2\bm{n}_{\mathfrak{u} -1}}) }\\
		& + \phi(f_{t_1 + \frac{1}{2n_1}})^{-A_1} \abs*{\Delta(\frac{f_{t_1 + \frac{1}{2n_1}}}{\phi(f_{t_1 + \frac{1}{2n_1}})^{-A_1}}; \bm{t}_{-1}, \bm{t}_{-1} + \frac{1}{2\bm{n}_{\mathfrak{u} -1}}) } d{t}_1 d\bm{t}_{-1}.
	\end{split}
\end{equation}
Notice that when $[{t}_1, {t}_1 + \frac{1}{2 {n}_{1}}] \cap \partial (\bar{\Omega}_{1}) \neq \varnothing$, we can verify the boundary growth condition of $\frac{f_{t_1 + \frac{1}{2n_1}}}{\phi(f_{t_1 + \frac{1}{2n_1}})^{-A_1}}$:
\begin{equation}
	\partial^{\mathfrak{u}} \frac{f_{t_1 + \frac{1}{2n_1}}}{\phi(f_{t_1 + \frac{1}{2n_1}})^{-A_1}} \lesssim \prod_{j = 2}^s \phi(t_j)^{-A_j - \mathbbm{1}_{j \in \mathfrak{u}}},
\end{equation}
for $\mathfrak{u} \subseteq 2:s$. The condition for $\frac{f_{t_1 + \frac{1}{2n_1}}}{\phi(f_{t_1 + \frac{1}{2n_1}})^{-A_1}}$ can be similarily verified. By our hypothesis, we have
\begin{equation}
	\int_{[0, 1]^{s-1}} \abs*{\Delta(\frac{f_{t_1}}{\phi(t_1)^{-A_1}}; \bm{t}_{-1}, \bm{t}_{-1} + \frac{1}{2\bm{n}_{\mathfrak{u} -1}}) } + \abs*{\Delta(\frac{f_{t_1 + \frac{1}{2n_1}}}{\phi(f_{t_1 + \frac{1}{2n_1}})^{-A_1}}; \bm{t}_{-1}, \bm{t}_{-1} + \frac{1}{2\bm{n}_{\mathfrak{u} -1}}) } \lesssim \prod_{j=2}^s {n}_j^{-1+A_j},
\end{equation}
when $[{t}_1, {t}_1 + \frac{1}{2 {n}_{1}}] \cap \partial (\bar{\Omega}_{1}) \neq \varnothing$. Next, we have
\begin{equation}
	\int_{0}^{1} \mathbbm{1}_{[{t}_1, {t}_1 + \frac{1}{2 {n}_{1}}] \cap \partial (\bar{\Omega}_{1}) \neq \varnothing} \phi(t_1)^{-A_1} \lesssim \int_{0}^{\frac{1}{2n_1}} t_1^{-A_1} dt_1 \lesssim n_1^{-1+A_1}.
\end{equation}
Thus,
\begin{equation}
	\int_{[0, 1]^{s-1}} \int_0^1 \mathbbm{1}_{[{t}_1, {t}_1 + \frac{1}{2 {n}_{1}}] \cap \partial (\bar{\Omega}_{1}) \neq \varnothing} \abs*{\Delta(f; \bm{t}, \bm{t} + \frac{1}{2\bm{n}_{\mathfrak{u}}}) } d {t}_1 d\bm{t}_{-1} \lesssim \prod_{j=1}^s {n}_j^{-1+A_j}.
\end{equation}
This concludes the first part of the proof. {Next, we consider the case when ${[{t}_1, {t}_1 + \frac{1}{2 {n}_{1}}] \cap \partial (\bar{\Omega}_{1}) = \varnothing}$. 
}
We have
\begin{equation}
	\begin{split}
		& \Delta(f_{t_1} - f_{t_1 + \frac{1}{2n_1}}; \bm{t}_{-1}, \bm{t}_{-1} + \frac{1}{2\bm{n}_{\mathfrak{u} -1}}) \\
		& = \Delta(\int_{[{t}_1, {t}_1 + \frac{1}{2 {n}_{1}} ]}\partial^{\tau} f_{\tau} d\tau; \bm{t}_{-1}, \bm{t}_{-1} + \frac{1}{2\bm{n}_{\mathfrak{u} -1}})\\
		&= \int_{t_1}^{t_1 + \frac{1}{2n_1}} \phi(\tau)^{-1-A_j} d\tau \Delta \left( \frac{\int_{[{t}_1, {t}_1 + \frac{1}{2 {n}_{1}} ]}\partial^{\tau} f_{\tau} d\tau}{\int_{t_1}^{t_1 + \frac{1}{2n_1}} \phi(\tau)^{-1-A_j} d\tau}; \bm{t}_{-1}, \bm{t}_{-1} + \frac{1}{2\bm{n}_{\mathfrak{u} -1}} \right).
	\end{split}
\end{equation}
It is straightforward to verify that $\frac{\int_{[{t}_1, {t}_1 + \frac{1}{2 {n}_{1}} ]}\partial^{\tau} f_{\tau} d\tau}{\int_{t_1}^{t_1 + \frac{1}{2n_1}} \phi(\tau)^{-1-A_j} d\tau}$ satisfies the boundary growth condition for all $ {[{t}_1, {t}_1 + \frac{1}{2 {n}_{1}}] \cap \partial (\bar{\Omega}_{1}) = \varnothing} $ in dimension $2:s$.

Thus, by our hypothesis, we have
\begin{equation}
	\int_{[0, 1]^{s-1}}  \abs*{\Delta\left(\frac{\int_{[{t}_1, {t}_1 + \frac{1}{2 {n}_{1}} ]}\partial^{\tau} f_{\tau} d\tau}{\int_{t_1}^{t_1 + \frac{1}{2n_1}} \phi(\tau)^{-1-A_j} d\tau}; \bm{t}_{-1}, \bm{t}_{-1} + \frac{1}{2\bm{n}_{\mathfrak{u} -1}}\right)} d\bm{t} \lesssim  \prod_{j=2}^s {n}_j^{-1+A_j}. 
\end{equation}
Finally, we obtain
\begin{equation}
	\label{eq:bound_cn_multidimension_axis_parallel_case_31}
	\begin{split}
	& \int_{[0, 1]^{s-1}} \int_0^1 \mathbbm{1}_{[{t}_1, {t}_1 + \frac{1}{2 {n}_{1}}] \cap \partial (\bar{\Omega}_{1}) = \varnothing} \abs*{\Delta(f; \bm{t}, \bm{t} + \frac{1}{2\bm{n}_{\mathfrak{u}}}) } d {t}_1 d\bm{t}_{-1}\\
	& = \int_{[0, 1]^{s-1}} \int_0^1 \mathbbm{1}_{[{t}_1, {t}_1 + \frac{1}{2 {n}_{1}}] \cap \partial (\bar{\Omega}_{1}) = \varnothing} \int_{t_1}^{t_1 + \frac{1}{2n_1}} \phi(\tau)^{-1-A_j} d\tau   \abs*{\Delta \left( \frac{\int_{[{t}_1, {t}_1 + \frac{1}{2 {n}_{1}} ]}\partial^{\tau} f_{\tau} d\tau}{\int_{t_1}^{t_1 + \frac{1}{2n_1}} \phi(\tau)^{-1-A_j} d\tau}; \bm{t}_{-1}, \bm{t}_{-1} + \frac{1}{2\bm{n}_{\mathfrak{u} -1}} \right)} d {t}_1 d\bm{t}_{-1}\\
	& \lesssim  \int_0^1  \mathbbm{1}_{[{t}_1, {t}_1 + \frac{1}{2 {n}_{1}}] \cap \partial (\bar{\Omega}_{1}) = \varnothing} \int_{[{t}_1, {t}_1 + \frac{1}{2{n}_1}  ]}  
	\phi(\tau)^{-1-A_1} d {\tau} d {t}_1 \prod_{j=2}^s {n}_j^{-1+A_j}  \lesssim \prod_{j=1}^s {n}_j^{-1+A_j}.
	\end{split}
\end{equation}

\end{proof}

\section{Walsh functions}
\label{section:walsh_functions}
Walsh (1923) introduced a system of functions denoted by ${}_{b}\wal_{k}$. 
\begin{Definition}[Walsh functions in 1-D]
	Let $k \in \mathbb{N}_0$ with $b$-adic expansion $k = \kappa_0 + \kappa_1 b + \kappa_2 b^2 + \cdots$. For $b \geq 2$ we denote by $\omega_b$ the primitive $b$th root of unity $\omega_b = e^{2\pi i/b}$. Then the $b$-ary Walsh function ${}_{b}\wal_{k}$ is defined by
	\begin{equation*}
		{}_{b}\wal_{k}(x) = \omega_b^{\kappa_0 \xi_1 + \kappa_1 \xi_2 + \cdots}
	\end{equation*}
	for $x \in [0, 1)$ with $b$-adic expansion $x = \xi_1 b^{-1} + \xi_2 b^{-2} + \cdots$ (unique in the sense that infinitely many of the digits $\xi_i$ must be different from $b-1$). 
\end{Definition}
The Walsh function in multiple dimensions is constructed by tensorizing the one-dimensional function.

In the following we present the proof of Lemma~\ref{lemma:walsh_series_T_ell}.
\begin{proof}
For $\bm{k} = ({k}_1,\dotsc, {k}_s)$, the $s$-dimensional Walsh-Dirichlet kernel is given by
\begin{equation}
	\label{eq:walsh_dirichlet_kernel}
	\begin{split}
		D_{\bm{\ell}}(\bm{z}) &= \sum_{{k}_1 = 0}^{b^{{\ell}_1} - 1} \cdots \sum_{{k}_s = 0}^{b^{{\ell}_s} - 1} {}_{b}\wal_{\bm{k}}(\bm{z}) = \sum_{{k}_1 = 0}^{b^{{\ell}_1} - 1} \cdots \sum_{{k}_s = 0}^{b^{{\ell}_s} - 1} \prod_{j=1}^s {}_{b}\wal_{{k}_j}({z}_j) = \prod_{j=1}^s b^{ {\ell}_j} \mathbbm{1}_{{z}_j \in [0, b^{- {\ell}_j})} = \mathbbm{1}_{\bm{z} \in [0, b^{-\bm{\ell}})} \prod_{j=1}^s b^{ {\ell}_j}. 
	\end{split}
\end{equation}
Notice that the Walsh--Fourier coefficient is given by
\begin{equation}
	\label{eq:walsh_fourier_coefficient}
	\bar{f}(\bm{k}) = \int_{[0, 1]^s} f(\bm{t}) \left({}_{b}\wal_{\bm{k}}(\bm{t})\right)^{*} d\bm{t}. 
\end{equation}
Substitute~\eqref{eq:walsh_fourier_coefficient} into~\eqref{eq:walsh_series_T_ell}, we have
\begin{equation}
	\begin{split}
		\sum_{\bm{k} \in T_{\bm{\ell}}} \bar{f}(\bm{k}) {}_{2} \wal_{\bm{k}}(\bm{t}) &= \sum_{\bm{k} \in T_{\bm{\ell}}} \int_{[0, 1]^s} f(\bm{y}) \left({}_{2}\wal_{\bm{k}}(\bm{y})\right)^{*} {}_{2} \wal_{\bm{k}}(\bm{t}) d\bm{y}\\
		&= \sum_{\bm{k} \in T_{\bm{\ell}}} \int_{[0, 1]^s} f(\bm{y}){}_{2}\wal_{\bm{k}}(\bm{y} \ominus \bm{t}) d\bm{y}\\
		&= \int_{[0, 1]^s} f(\bm{y}) \sum_{\bm{k} \in T_{\bm{\ell}}} {}_{2}\wal_{\bm{k}}(\bm{y} \ominus \bm{t}) d\bm{y},
	\end{split}
\end{equation}
where $\ominus$ denotes the digit-wise subtraction modulo $2$. Refering to~\eqref{eq:walsh_dirichlet_kernel}, the Walsh--Dirichlet kernel $\sum_{\bm{k} \in T_{\bm{\ell}}} {}_{2}\wal_{\bm{k}}(\bm{y} \ominus \bm{t})$ is nonzero  when $\bm{y} \ominus \bm{t} \in [0, 2^{-\bm{\ell}}]$. In this case, the first ${\ell}_j$ digits of ${y}_j$ and ${t}_j$ must match for all $j = 1,\dotsc, s$, which is equivalent to the condition $\left\lfloor {y}_j 2^{{\ell}_j} \right\rfloor = \left\lfloor {t}_j 2^{{\ell}_j} \right\rfloor$ for $j = 1, \dotsc, s$. Thus we have
\begin{equation}
	\sum_{\bm{k} \in T_{\bm{\ell}}} \bar{f}(\bm{k}) {}_{2} \wal_{\bm{k}}(\bm{t}) = \prod_{j=1}^s 2^{{\ell}_j} \int_{\cap_{j=1}^s \left\lfloor {y}_j 2^{{\ell}_j} \right\rfloor = \left\lfloor {t}_j 2^{{\ell}_j} \right\rfloor} f(\bm{y}) d\bm{y}. 
\end{equation}
\end{proof}

\section{Proof of Lemma~\ref{lemma:sobol_axis_parallel_multiple}}
\label{sec:proof_sobol_axis_parallel_multiple}
\begin{proof}
	Following~\cite{he2015convergence}, we consider a partition of $[0, 1]^s$ with the elementary intervals $E_{\bm{\ell}, \bm{k}}$, where $b=2$, $\bm{\ell} \in \mathbb{N}_0^s$, as defined in~\eqref{eq:elementary_interval}. Using notations similar to that in Lemma~\ref{lemma:alternating_sum_integration}, we define $\bar{\bm{\ell}} \coloneqq \max(\bm{\ell} - 1, 0)$, $\bar{\bm{\ell}}_{\mathfrak{u}} \coloneqq {\bm{\ell}}_{\mathfrak{u}} - 1$ and $f(\bm{y}) \coloneqq 0$ when $\bm{y} \notin [0, 1]^s$ for simplicity. For a given $\bm{\ell} \in \mathbb{N}_0^s$, we have
\begin{equation}
	\sigma_{\bm{\ell}}^2(f_M) = \sum_{\substack{\bm{k}  \in \mathbb{N}_0^s  \\ {\bm{k} \leq 2^{\bar{\bm{\ell}}} - 1} }} \left(\int_{E_{\bm{\ell}, 2\bm{k}}} {\Delta(f_M; \bm{y}, \bm{y} + \frac{1}{2^{\bm{\ell}} })} d\bm{y} \right)^2.  
	\label{eq:sigma2_ell_multi_dimension_fM}
\end{equation}
In the following we expand~\eqref{eq:sigma2_ell_multi_dimension_fM} by considering all possible cases where the elementary interval $E_{\bar{\bm{\ell}}, \bm{k}}$ intersects with the boundary of $\Omega$. Specifically, we have
\begin{equation}
	\sigma_{\bm{\ell}}^2(f_M) = \sum_{\mathfrak{u} \subseteq 1:s} \sum_{\substack{\bm{k}  \in \mathbb{N}_0^s  \\ {\bm{k} \leq 2^{\bar{\bm{\ell}}} - 1} }} \mathbbm{1}_{\substack{\cap_{j \in \mathfrak{u}} \left( (E_{\bar{\bm{\ell}}, \bm{k}})_j \cap  \partial (\Omega_j) \neq \varnothing \right) \\ \cap_{i \notin \mathfrak{u} } \left( (E_{\bar{\bm{\ell}}, \bm{k}})_i \cap \partial (\Omega_i)  = \varnothing \right) } } \left(\int_{E_{\bm{\ell}, 2\bm{k}}} {\Delta(f_M; \bm{y}, \bm{y} + \frac{1}{2^{\bm{\ell}} })} d\bm{y} \right)^2,  
\end{equation}
where $(E_{\bar{\bm{\ell}}, \bm{k}})_j$ denotes the projection of the interval $(E_{\bar{\bm{\ell}}, \bm{k}})$ onto dimension j, $\partial (\Omega_j)$ denotes the boundary of the projection of $\Omega$ onto the dimension $j$, and $\mathbbm{1}_{\substack{\cap_{j \in \mathfrak{u}} \left( (E_{\bar{\bm{\ell}}, \bm{k}})_j \cap \partial (\Omega_j) \neq \varnothing \right) \\ \cap_{i \notin \mathfrak{u} } \left( (E_{\bar{\bm{\ell}}, \bm{k}})_i \cap \partial (\Omega_i)  = \varnothing \right) } }$ takes value 1 if the event $\cap_{j \in \mathfrak{u}} \left( (E_{\bar{\bm{\ell}}, \bm{k}})_j \cap \partial (\Omega_j) \neq \varnothing \right) \cap_{i \notin \mathfrak{u} } \left( (E_{\bar{\bm{\ell}}, \bm{k}})_i \cap \partial (\Omega_i)  = \varnothing \right)$ is true and 0 otherwise. 
For a given $\mathfrak{u} \subseteq 1:s$, we have
\begin{equation}
	\begin{split}
		& \sum_{\substack{\bm{k}  \in \mathbb{N}_0^s  \\ {\bm{k} \leq 2^{\bar{\bm{\ell}}} - 1} }} \mathbbm{1}_{\substack{\cap_{j \in \mathfrak{u}} \left( (E_{\bar{\bm{\ell}}, \bm{k}})_j \cap \partial (\Omega_j) \neq \varnothing \right) \\ \cap_{i \notin \mathfrak{u} } \left( (E_{\bar{\bm{\ell}}, \bm{k}})_i \cap \partial (\Omega_i)  = \varnothing \right) } } \left(\int_{E_{\bm{\ell}, 2\bm{k}}} {\Delta(f_M; \bm{y}, \bm{y} + \frac{1}{2^{\bm{\ell}} })} d\bm{y} \right)^2 \\
		&=  \sum_{\substack{\bm{k}_{\mathfrak{u}}  \in \mathbb{N}_0^{\abs{\mathfrak{u}}}  \\ {\bm{k}_{\mathfrak{u}} \leq 2^{\bar{\bm{\ell}}_{\mathfrak{u}} } - 1} }} \mathbbm{1}_{{\cap_{j \in \mathfrak{u}} \left( (E_{\bar{\bm{\ell}}, \bm{k}})_j \cap \partial (\Omega_j) \neq \varnothing \right) } } \sum_{\substack{\bm{k}_{-\mathfrak{u}}  \in \mathbb{N}_0^{\abs{-\mathfrak{u}}}  \\ {\bm{k}_{-\mathfrak{u}} \leq 2^{\bar{\bm{\ell}}_{-\mathfrak{u}} } - 1} }}  \mathbbm{1}_{\cap_{i \notin \mathfrak{u} } \left( (E_{\bar{\bm{\ell}}, \bm{k}})_i \cap \partial (\Omega_i)  = \varnothing \right) }  \left(\int_{E_{\bm{\ell}, 2\bm{k}}} {\Delta(f_M; \bm{y}, \bm{y} + \frac{1}{2^{\bm{\ell}} })} d\bm{y} \right)^2.
	\end{split}
\end{equation}
For $\mathfrak{u} \subseteq 1:s$, we have the following decomposition:
\begin{equation}
		{\Delta(f_M; \bm{y}, \bm{y} + \frac{1}{2^{\bm{\ell}} })} = \sum_{\mathfrak{v} \subseteq \mathfrak{u}} (-1)^{\abs{\mathfrak{v}}} \sum_{\mathfrak{z} \subseteq -\mathfrak{u}} (-1)^{\abs{\mathfrak{z}}} f_M ( \bm{y}^{\mathfrak{z}}; (\bm{y} + 2^{-\bm{\ell}})^{-\mathfrak{u}-\mathfrak{z}} \mid \bm{y}^{\mathfrak{v}}; (\bm{y} + 2^{-\bm{\ell}})^{\mathfrak{u}-\mathfrak{v}} ),
\end{equation}
where $ \bm{y}^{\mathfrak{z}}; (\bm{y} + 2^{-\bm{\ell}})^{-\mathfrak{u}-\mathfrak{z}} \mid \bm{y}^{\mathfrak{v}}; (\bm{y} + 2^{-\bm{\ell}})^{\mathfrak{u}-\mathfrak{v}} $ denotes a $\bm{y} \in [0, 1]^s$ such that $\bm{y}_{\mathfrak{u}} = \bm{y}^{\mathfrak{v}}; (\bm{y} + 2^{-\bm{\ell}})^{\mathfrak{u}-\mathfrak{v}}$ and $\bm{y}_{-\mathfrak{u}} = \bm{y}^{\mathfrak{z}}; (\bm{y} + 2^{-\bm{\ell}})^{-\mathfrak{u}-\mathfrak{z}}$. 
We have the following triangular inequality,
\begin{equation}
	\abs*{\Delta(f_M; \bm{y}, \bm{y} + \frac{1}{2^{\bm{\ell}} })} \leq \sum_{\mathfrak{v} \subseteq \mathfrak{u}}  \abs*{\sum_{\mathfrak{z} \subseteq -\mathfrak{u}} (-1)^{\abs{\mathfrak{z}}} f_M ( \bm{y}^{\mathfrak{z}}; (\bm{y} + 2^{-\bm{\ell}})^{-\mathfrak{u}-\mathfrak{z}} \mid \bm{y}^{\mathfrak{v}}; (\bm{y} + 2^{-\bm{\ell}})^{\mathfrak{u}-\mathfrak{v}} )},
\end{equation}
and thus
\begin{equation} 
	\begin{split}
	& \left(\int_{E_{\bm{\ell}, 2\bm{k}}} {\Delta(f_M; \bm{y}, \bm{y} + \frac{1}{2^{\bm{\ell}} })} d\bm{y} \right)^2 
	\leq \left(\int_{E_{\bm{\ell}, 2\bm{k}}} \abs*{\Delta(f_M; \bm{y}, \bm{y} + \frac{1}{2^{\bm{\ell}} })} d\bm{y} \right)^2 \\
	&\leq \left(\int_{E_{\bm{\ell}, 2\bm{k}}} \sum_{\mathfrak{v} \subseteq \mathfrak{u}}  \abs*{\sum_{\mathfrak{z} \subseteq -\mathfrak{u}} (-1)^{\abs{\mathfrak{z}}} f_M ( \bm{y}^{\mathfrak{z}}; (\bm{y} + 2^{-\bm{\ell}})^{-\mathfrak{u}-\mathfrak{z}} \mid \bm{y}^{\mathfrak{v}}; (\bm{y} + 2^{-\bm{\ell}})^{\mathfrak{u}-\mathfrak{v}} )}  d\bm{y} \right)^2\\
	&\leq 2^{\abs{\mathfrak{u}}} \sum_{\mathfrak{v} \subseteq \mathfrak{u}} \left(\int_{E_{\bm{\ell}, 2\bm{k}}}   \abs*{\sum_{\mathfrak{z} \subseteq -\mathfrak{u}} (-1)^{\abs{\mathfrak{z}}} f_M ( \bm{y}^{\mathfrak{z}}; (\bm{y} + 2^{-\bm{\ell}})^{-\mathfrak{u}-\mathfrak{z}} \mid \bm{y}^{\mathfrak{v}}; (\bm{y} + 2^{-\bm{\ell}})^{\mathfrak{u}-\mathfrak{v}} )}  d\bm{y} \right)^2.
	\end{split}
\end{equation}
We can split the integration domain into the dimensions $\mathfrak{u}$ and $-\mathfrak{u}$ using the following inequality:
\begin{equation}
	\begin{split}
		 &  \int_{E_{\bm{\ell}, 2\bm{k}}}   \abs*{\sum_{\mathfrak{z} \subseteq -\mathfrak{u}} (-1)^{\abs{\mathfrak{z}}} f_M ( \bm{y}^{\mathfrak{z}}; (\bm{y} + 2^{-\bm{\ell}})^{-\mathfrak{u}-\mathfrak{z}} \mid \bm{y}^{\mathfrak{v}}; (\bm{y} + 2^{-\bm{\ell}})^{\mathfrak{u}-\mathfrak{v}} )}  d\bm{y} \leq \int_{(E_{\bm{\ell}, 2\bm{k}})_{\mathfrak{u}} } \phi (\bm{y}^{\mathfrak{v}}; (\bm{y} + 2^{-\bm{\ell}})^{\mathfrak{u}-\mathfrak{v}})^{-A^{*}}  d\bm{y}_{ \mathfrak{u}} \\
		 &\cdot \sup_{\bm{y}_{\mathfrak{u}} \in (E_{\bm{\ell}, 2\bm{k}})_{\mathfrak{u}}} \int_{(E_{\bm{\ell}, 2\bm{k}})_{-\mathfrak{u}} } \frac{\abs*{\Delta \left( f_M ( \cdot \mid \bm{y}^{\mathfrak{v}}; (\bm{y} + 2^{-\bm{\ell}})^{\mathfrak{u}-\mathfrak{v}} ); \bm{y}, \bm{y} + 2^{-\bm{\ell}_{-\mathfrak{u}}} \right) }}{\phi (\bm{y}^{\mathfrak{v}}; (\bm{y} + 2^{-\bm{\ell}})^{\mathfrak{u}-\mathfrak{v}})^{-A^{*}} } d\bm{y}_{-\mathfrak{u}}. 
	\end{split}
\end{equation}
The boundary growth condition of $\frac{{   f_M ( \cdot \mid \bm{y}^{\mathfrak{v}}; (\bm{y} + 2^{-\bm{\ell}})^{\mathfrak{u}-\mathfrak{v}} ) }}{\phi (\bm{y}^{\mathfrak{v}}; (\bm{y} + 2^{-\bm{\ell}})^{\mathfrak{u}-\mathfrak{v}})^{-A^{*}} }$ can be verified and from Lemma~\ref{lemma:variance_sobol_multi_d}, we have
\begin{equation}
	\begin{split}
		& \sum_{\substack{\bm{k}_{-\mathfrak{u}}  \in \mathbb{N}_0^{\abs{-\mathfrak{u}}}  \\ {\bm{k}_{-\mathfrak{u}} \leq 2^{\bar{\bm{\ell}}_{-\mathfrak{u}} } - 1} }}  \mathbbm{1}_{\cap_{i \notin \mathfrak{u} } \left( (E_{\bar{\bm{\ell}}, \bm{k}})_i \cap \partial (\Omega)_i  = \varnothing \right) }\\
		& \left( \sup_{\bm{y}_{\mathfrak{u}}\in (E_{\bm{\ell}, 2\bm{k}})_{\mathfrak{u}}} \int_{(E_{\bm{\ell}, 2\bm{k}})_{-\mathfrak{u}} } \frac{\abs*{\Delta \left( f_M ( \cdot \mid \bm{y}^{\mathfrak{v}}; (\bm{y} + 2^{-\bm{\ell}})^{\mathfrak{u}-\mathfrak{v}} ); \bm{y}, \bm{y} + 2^{-\bm{\ell}_{-\mathfrak{u}}} \right) }}{\phi (\bm{y}^{\mathfrak{v}}; (\bm{y} + 2^{-\bm{\ell}})^{\mathfrak{u}-\mathfrak{v}})^{-A^{*}} } d\bm{y}_{-\mathfrak{u}} \right)^2 \lesssim 2^{(2A^{*} - 1)\abs*{\bm{\ell}_{-\mathfrak{u}}}}.
	\end{split}
\end{equation}
Moreover, for all $\bm{k}_{-\mathfrak{u}}  \in \mathbb{N}_0^{\abs{-\mathfrak{u}}}$, we have
\begin{equation}
	\sum_{\substack{\bm{k}_{\mathfrak{u}}  \in \mathbb{N}_0^{\abs{\mathfrak{u}}}  \\ {\bm{k}_{\mathfrak{u}} \leq 2^{\bar{\bm{\ell}}_{\mathfrak{u}} } - 1} }} \mathbbm{1}_{{\cap_{j \in \mathfrak{u}} \left( (E_{\bar{\bm{\ell}}, \bm{k}})_j \cap \partial (\Omega)_j \neq \varnothing \right) } } \sum_{\mathfrak{v} \subseteq \mathfrak{u}} \left( \int_{(E_{\bm{\ell}, 2\bm{k}})_{\mathfrak{u}} } \phi (\bm{y}^{\mathfrak{v}}; (\bm{y} + 2^{-\bm{\ell}})^{\mathfrak{u}-\mathfrak{v}})^{-A^{*}}  d\bm{y}_{ \mathfrak{u}} \right)^2 \lesssim M 2^{\abs{\bm{\ell}_{\mathfrak{u}} }}. 
\end{equation}
Thus we have,
\begin{equation}
	\begin{split}
		\sigma^2_{\bm{\ell}} (f_M) &\leq \sum_{\mathfrak{u}\subseteq 1:s} 2^{\abs{\mathfrak{u}}} \sum_{\substack{\bm{k}_{\mathfrak{u}}  \in \mathbb{N}_0^{\abs{\mathfrak{u}}}  \\ {\bm{k}_{\mathfrak{u}} \leq 2^{\bar{\bm{\ell}}_{\mathfrak{u}} } - 1} }} \mathbbm{1}_{{\cap_{j \in \mathfrak{u}} \left( (E_{\bar{\bm{\ell}}, \bm{k}})_j \cap \partial (\Omega)_j \neq \varnothing \right) } } \left( \int_{(E_{\bm{\ell}, 2\bm{k}})_{\mathfrak{u}} } \phi (\bm{y}^{\mathfrak{v}}; (\bm{y} + 2^{-\bm{\ell}})^{\mathfrak{u}-\mathfrak{v}})^{-A^{*}}  d\bm{y}_{ \mathfrak{u}} \right)^2   \\
		& \sum_{\substack{\bm{k}_{-\mathfrak{u}}  \in \mathbb{N}_0^{\abs{-\mathfrak{u}}}  \\ {\bm{k}_{-\mathfrak{u}} \leq 2^{\bar{\bm{\ell}}_{-\mathfrak{u}} } - 1} }} \mathbbm{1}_{\cap_{i \notin \mathfrak{u} } \left( (E_{\bar{\bm{\ell}}, \bm{k}})_i \cap \partial (\Omega)_i  = \varnothing \right) } \sum_{\mathfrak{v} \subseteq \mathfrak{u}} \left( \sup_{\bm{y}_{\mathfrak{u}}} \int_{(E_{\bm{\ell}, 2\bm{k}})_{-\mathfrak{u}} } \frac{\abs*{\Delta \left( f_M ( \cdot \mid \bm{y}^{\mathfrak{v}}; (\bm{y} + 2^{-\bm{\ell}})^{\mathfrak{u}-\mathfrak{v}} ); \bm{y}, \bm{y} + 2^{-\bm{\ell}_{-\mathfrak{u}}} \right) }}{\phi (\bm{y}^{\mathfrak{v}}; (\bm{y} + 2^{-\bm{\ell}})^{\mathfrak{u}-\mathfrak{v}})^{-A^{*}} } d\bm{y}_{-\mathfrak{u}} \right)^2 \\
		& \lesssim  M \sum_{\mathfrak{u}\subseteq 1:s} 2^{\abs{\mathfrak{u}}} 2^{-\abs{\bm{\ell}_{\mathfrak{u}}} } 2^{(2A^{*} - 1)\abs*{\bm{\ell}_{-\mathfrak{u}}} } \lesssim M2^{\max(-1, 2A^{*} - 1)\abs{\bm{\ell}}}.
	\end{split}
\end{equation}
This proves Lemma~\ref{lemma:sobol_axis_parallel} when $M = 1$. 

Finally, we have,
\begin{equation}
	\var{I_N^{\mathrm{sob}}(f_M) } \lesssim M \left( N^{2A^{*} - 2} (\log N)^{s-1} +  N^{-2} (\log N)^{s-1} \right).
\end{equation}
\end{proof}

\section{Proof of Lemma~\ref{lemma:sobol_partial_axis_parallel_variance_decay} }
\label{section:proof_sobol_partial_axis_parallel}
\begin{proof}
Given the definition $\Omega = \Omega_{\mathfrak{u}} \times [\bm{a}_{-\mathfrak{u}}, \bm{b}_{-\mathfrak{u}}]$, $\mathfrak{u} \subsetneq 1:s$,
for $\bm{r}$ and $r$ defined in~\eqref{eq:elementary_interval_r}, we define ${\mathcal{T}_{\mathrm{tot}}^{r_{\mathfrak{u}}} } \coloneqq \{ \bigcup E_{\bm{r}, \bm{k}_{\mathfrak{u}}} : E_{\bm{r}, \bm{k}_{\mathfrak{u}}} \cap \Omega_{\mathfrak{u}} \neq \varnothing, 0 \leq {k}_j \leq {r} - 1, j \in \mathfrak{u} \}$. 
Define $f_{r_{\mathfrak{u}}} \coloneqq \mathbbm{1}_{{\mathcal{T}_{\mathrm{tot}}^{r_{\mathfrak{u}}} } \times [\bm{a}_{-\mathfrak{u}}, \bm{b}_{-\mathfrak{u}}]} g$. We have the following inequality
	\begin{equation}
		\var{I_N^{\mathrm{sob}} (f)} \leq 2 \var{I_N^{\mathrm{sob}}(f - f_{{r}_{\mathfrak{u}}})} + 2 \var{I_N^{\mathrm{sob}}(f_{{r}_{\mathfrak{u}}})}. 
	 \end{equation}
	 Following Lemma~\ref{lemma:sobol_axis_parallel_multiple} and the derivations in Section~\ref{sec:general_discontinuous_integrand_sobol}, we have
	 \begin{equation}
		\var{I_N^{\mathrm{sob}} (f_{{r}_{\mathfrak{u}}})} \lesssim r^{\abs{\mathfrak{u}} - 1} \left( N^{-2} (\log N)^{s - 1} + N^{2A^{*} - 2} (\log N)^{s - 1} \right).
	 \end{equation}
	 When $A^{*} < 0$, we have
	 \begin{equation}
		\var{I_N^{\mathrm{sob}}(f - f_{{r}_{\mathfrak{u}}})} \leq N^{-1} \E{(f - f_{{r}_{\mathfrak{u}}})^2} \lesssim N^{-1}  r^{\abs{\mathfrak{u}} - 1} r^{-\abs{\mathfrak{u}}} = N^{-1} r^{-1},
	 \end{equation}
	and
	\begin{equation}
		\var{I_N^{\mathrm{sob}} (f)} \lesssim N^{-1} r^{-1} + r^{\abs{\mathfrak{u}} - 1} N^{-2} (\log N)^{s - 1}.
	\end{equation}
	We choose $r = N^{\frac{1}{ \abs{\mathfrak{u}} }} (\log N)^{\frac{1-s}{ \abs{\mathfrak{u}} }}$ to obtain
	\begin{equation}
		\var{I_N^{\mathrm{sob}} (f)} \lesssim N^{-1 - \frac{1}{ \abs{\mathfrak{u}}  }} (\log N)^{\frac{s - 1}{ \abs{\mathfrak{u}} }}. 
	\end{equation}
	When $A^{*} \geq 0$, we have
	 \begin{equation}
		\begin{split}
			\var{I_N^{\mathrm{sob}}(f - f_{{r}_{\mathfrak{u}}})} \leq N^{-1} \E{(f - f_{r_{\mathfrak{u}}})^2} &\leq N^{-1} r^{-1 + 2A^{*} + \delta},
		\end{split}
	\end{equation}
	for $0 < \delta < 1-2A^{*}$. We have
	\begin{equation}
		\var{I_N^{\mathrm{sob}} (f)} \lesssim N^{-1} r^{-1 + 2A^{*} + \delta} + r^{ \abs{\mathfrak{u}} - 1}  N^{2A^{*} - 2} (\log N)^{s - 1}.
	\end{equation}
	We choose $r = N^{-\frac{1-2A^{*}}{ 2A^{*} + \delta -  \abs{\mathfrak{u}}}} (\log N)^{\frac{s-1}{ 2A^{*} + \delta -  \abs{\mathfrak{u}}}}$ to obtain
	\begin{equation}
		\var{I_N^{\mathrm{sob}} (f)} \lesssim N^{-1 + \frac{(1-2A^{*})(-1 + 2A^{*} + \delta)}{ \abs{\mathfrak{u}}  - 2A^{*} - \delta}} (\log N)^{ \frac{(s-1)(1 - 2A^{*} - \delta) }{ \abs{\mathfrak{u}}  - 2A^{*} - \delta}}.
	\end{equation}
\end{proof}